\documentclass[a4paper,10pt,reqno]{amsart}
\usepackage{graphics}
\newtheorem{lem}{Lemma}

\newtheorem{theo}{Theorem}

\newtheorem{cor}{Corollary}
\newtheorem{defi}{Definition}
\newtheorem{prop}{Proposition}
\numberwithin{equation}{section}

\newcommand{\sgn}{\operatorname{sgn}}

\newcommand{\id}{\operatorname{id}}

\begin{document}

\title[Gelfand-Tsetlin tree sequences]
{Sequences of labeled trees related to Gelfand-Tsetlin patterns}

\author[Ilse Fischer]{Ilse Fischer}

\thanks{Supported by the Austrian Science Foundation
    FWF, START grant Y463 and NFN grant S9607--N13.}

\begin{abstract}
By rewriting the famous hook-content formula it easily follows that there are $\prod\limits_{1 \le i < j \le n} \frac{k_j - k_i + j - i}{j-i}$ semistandard tableaux of shape $(k_n,k_{n-1},\ldots,k_1)$ with entries in $\{1,2,\ldots,n\}$ or, equivalently, Gelfand-Tsetlin patterns with bottom row $(k_1,\ldots,k_n)$. In this article we introduce certain sequences of labeled trees, the signed enumeration of which is also given by this formula. In these trees, vertices as well as edges are labeled, the crucial condition being that each edge label lies between the vertex labels of the two endpoints of the edge. This notion enables us to give combinatorial explanations of the shifted antisymmetry of the formula 
and its polynomiality. Furthermore, we propose to develop an analog approach of combinatorial reasoning for monotone triangles  and explain how this may lead to a combinatorial understanding of the alternating sign matrix theorem.
\end{abstract}

\maketitle

\section{Introduction}
\label{intro}

One possibility to see that the expression
\begin{equation}
\label{product}
\prod\limits_{1 \le i < j \le n} \frac{k_j - k_i + j - i}{j-i}
\end{equation}  
is an integer for any choice of $(k_1,\ldots,k_n) \in \mathbb{Z}^n$ 
is to find combinatorial 
objects that are enumerated by this quantity. This is, for instance, accomplished by {\it Gelfand-Tsetlin patterns} 
with prescribed bottom row $k_1,k_2,\ldots,k_n$. A Gelfand-Tsetlin pattern (see \cite[p. 313]{stan} or \cite[(3)]{gel}  for the first appearance) is a triangular array of 
integers with $n$ rows of the following shape 
\begin{center}
\begin{tabular}{ccccccccccccccccc}
  &   &   &   &   &   &   &   & $a_{1,1}$ &   &   &   &   &   &   &   & \\
  &   &   &   &   &   &   & $a_{2,1}$ &   & $a_{2,2}$ &   &   &   &   &   &   & \\
  &   &   &   &   &   & $\dots$ &   & $\dots$ &   & $\dots$ &   &   &   &   &   & \\
  &   &   &   &   & $a_{n-2,1}$ &   & $\dots$ &   & $\dots$ &   & $a_{n-2,n-2}$ &   &   &   &   & \\
  &   &   &   & $a_{n-1,1}$ &   & $a_{n-1,2}$ &  &   $\dots$ &   & $\dots$   &  & $a_{n-1,n-1}$  &   &   &   & \\
  &   &   & $a_{n,1}$ &   & $a_{n,2}$ &   & $a_{n,3}$ &   & $\dots$ &   & $\dots$ &   & $a_{n,n}$ &   &   &
\end{tabular},
\end{center}
that is monotone increasing along northeast diagonals and southeast
diagonals, i.e. $a_{i,j} \le
a_{i-1,j}$ for $1 \le j  <  i \le n$ and  $a_{i,j} \le a_{i+1,j+1}$ for
$1 \le  j  \le i < n$. It is conceivable to assume that
$(k_1,\ldots,k_n) \in \mathbb{Z}_{\ge 0}^n$, as Gelfand-Tsetlin patterns with bottom row $(k_1,\ldots,k_n)$ are obviously in bijective correspondence with Gelfand-Tsetlin patterns with bottom row $(k_1+t,\ldots,k_n+t)$ for any integer $t \in \mathbb{Z}$. Under this assumption, they are equivalent to
{\it semistandard tableaux} of shape $(k_n,k_{n-1},\ldots,k_1)$ with entries in $\{1,2,\ldots,n\}$, the latter
being fillings of the Ferrers diagram associated with the integer partition $(k_n,k_{n-1},\ldots,k_1)$ that are weakly increasing 
	along rows and strictly increasing along columns.\footnote{Note that there is actually no dependency between the number of  feasible values for the entries of the semistandard tableaux and the number of parts in the integer partition: semistandard tableaux of shape 
		$(k_m,k_{m-1},\ldots,k_1)$ with entries in $\{1,2,\ldots,n\}$  are equivalent to semistandard tableaux of shape $(k_m,k_{m-1},\ldots,k_1,0^{n-m})$ with entries in $\{1,2,\ldots,n\}$ if $n \ge m$ and there exists no semistandard tableau otherwise.}  Next we give an example of a Gelfand-Tsetlin pattern and 
the corresponding semistandard tableaux.
\begin{center}
\begin{tabular}{ccccccccccc}
   &   &   &   &   & $2$ &   &   &   &   &   \\
   &   &   &   & $2$ &   & $2$ &   &   &   &    \\
   &   &   & $1$ &   & $2$ &   & $4$ &   &   &   \\
   &   & $1$ &   & $1$ &   & $3$ &   & $4$ &   &   \\
   & $0$ &   & $1$ &  &   $3$ &   & $3$   &  & $5$  &    \\
   $0$ &   & $0$ &   & $2$ &   & $3$ &   & $5$ &   & $6$ 
\end{tabular}
\qquad
\begin{tabular}{cccccc}
1 & 1 & 3 & 3 & 5 & 6 \\
2 & 2 & 4 & 6 & 6 &    \\
3 & 5 & 5 &    &    &    \\
4 & 6 &    &    &    &    
\end{tabular}
\end{center}
In general, given a Gelfand-Tsetlin pattern $(a_{i,j})_{1 \le j \le i \le n}$, the corresponding semistandard tableau is constructed by placing the integer $i$ in the cells of the skew shape
$$(a_{i,i},a_{i,i-1},\ldots,a_{i,1}) / (a_{i-1,i-1},a_{i-1,i-2},\ldots,a_{i-1,1}).$$

\medskip

Semistandard tableaux of fixed shape (and thus Gelfand-Tsetlin patterns) are known to be enumerated by the hook-content formula \cite[Corollary 7.21.4]{stan}, which is easily seen to be equivalent to \eqref{product}, see also \cite[Lemma 7.21.1]{stan}. A common way to prove this formula is to translate the problem into the enumeration of families of non-intersecting lattice paths with a certain set of fixed starting points and end points. To complement the treatment given in this article, we sketch this point of view in Appendix~\ref{path}.
A direct proof of the fact that Gelfand-Tsetlin patterns with bottom row 
$k_1, k_2, \ldots, k_n$ are enumerated by \eqref{product} can be found in \cite[Section~5]{method}. There we have actually proven a more general result, which 
we describe in the following paragraph.

\medskip

The reader will have noticed that the combinatorial interpretations that we have given so far only provide an explanation for the integrality of \eqref{product} if the sequence $k_1,k_2,\ldots,k_n$ is weakly increasing.  This can be overcome\footnote{Of course, this also follows by choosing a permutation 
$\sigma \in {\mathcal S}_n$ with $k_{\sigma_1}  + \sigma_1 \le k_{\sigma_2} + \sigma_2 \le \ldots \le k_{\sigma_n} + \sigma_n$ and then observing that 
$\prod\limits_{1 \le i < j \le n} \frac{k_{\sigma_j} - k_{\sigma_i} + \sigma_j - \sigma_i}{j-i}  = \sgn \sigma 
\prod\limits_{1 \le i < j \le n} \frac{k_j - k_i + j - i}{j-i}$ is the number of Gelfand-Tsetlin patterns with bottom row $(k_{\sigma_1}+\sigma_1-1,k_{\sigma_2}+\sigma_2-2,\ldots,k_{\sigma_n}+\sigma_n-n)$.} by extending the combinatorial interpretation of Gelfand-Tsetlin patterns with bottom row $(k_1,\dots,k_n)$ to all $n$-tuples of integers 
$(k_1,\ldots,k_n)$ and working with a signed enumeration as follows: a (generalized) Gelfand-Tsetlin pattern is an array of integers $(a_{i,j})_{1 \le j \le i \le n}$ such that the following condition is fulfilled: for any $a_{i,j}$ with $1 \le j \le i \le n-1$ we have $a_{i+1,j} \le a_{i,j} \le a_{i+1,j+1}$ if $a_{i+1,j} \le a_{i+1,j+1}$ and $a_{i+1,j} > a_{i,j} > a_{i+1,j+1}$
if $a_{i+1,j} > a_{i+1,j+1}$. (In particular, there exists no generalized Gelfand-Tsetlin pattern with $a_{i+1,j} = a_{i+1,j+1}+1$.) In the latter case we say that $a_{i,j}$ is an inversion. 
The weight (or sign) of a given Gelfand-Tsetlin pattern is $(-1)^\text{$\#$ of inversions}$. With this, \eqref{product} is the signed enumeration of all Gelfand-Tsetlin patterns with bottom row $k_1, k_2, \ldots, k_n$.

\medskip

The main task of the present paper is to provide a whole family of sets of objects that come along with a rather canonical notion of a sign, the signed enumeration of each of these sets is given by \eqref{product}. We call these objects {\it Gelfand-Tsetlin tree sequences} as Gelfand-Tsetlin patterns are one special member of this family. The definition of these objects is given in Section~\ref{definition}.
This enables us to give a combinatorial proof  of the fact that Gelfand-Tsetlin patterns are enumerated by \eqref{product}. Interestingly, this combinatorial 
proof is not based on a bijection between Gelfand-Tsetlin patterns and a second type of objects which are more easily seen to be enumerated by \eqref{product}. Rather than that we give combinatorial proofs of the facts that the replacement $(k_{i},k_{j}) \to (k_{j}+j-i,k_i-i+j)$ in the enumeration formula 
for the number of Gelfand-Tsetlin patterns with prescribed bottom row only causes the inversion of the sign (Section~\ref{indsym}) as well as that the enumeration formula must be a polynomial in $(k_1,\ldots,k_n)$ of degree no greater than $n-1$ in every $k_i$ (Section~\ref{diff}). For each of these properties, this is accomplished by providing an appropriate member of the family for which the respective property is almost obvious. Then, it is not hard to see that these properties essentially determine the enumeration formula, which is the only algebraic part of the proof.
Note that the first property can obviously only be understood combinatorially after having extended the combinatorial interpretation of Gelfand-Tsetlin patterns with bottom row $k_1,k_2,\ldots,k_n$ to 
arbitrary $(k_1,k_2,\ldots,k_n) \in \mathbb{Z}^n$ as the sequence
$k_1,\ldots,k_{i-1},k_j+j-i,k_{i+1},\ldots,k_{j-1},k_i+i-j,k_{j+1},\ldots,k_n$ can not be weakly increasing if $k_1,k_2,\ldots,k_n$ is weakly increasing. Also the inversion of the sign surely indicates that a signed enumeration must be involved.

\medskip

However, the original motivation for this paper is the intention to 
translate some of the research we have done on {\it monotone triangles} into a more combinatorial reasoning. Monotone triangles are Gelfand-Tsetlin patterns with strictly increasing rows and their significance is due to the fact that they are in bijective correspondence with alternating sign matrices when prescribing $1,2,\ldots,n$ as bottom row.
It took a lot of effort to enumerate $n \times n$ alternating sign matrices and all proofs known so far can not be considered as combinatorial proofs as they usually involve heavy algebraic manipulations, see \cite{bressoud}. Also the long-standing ``Gog-Magog conjecture'' \cite{kratt}, which is a generalization of the fact that $n \times n$ alternating sign matrices are in bijective correspondence with $2n \times 2n \times 2n$ totally symmetric self-complementary plane partitions is still unsolved, which is another indication for the fact that alternating sign matrices (as well as plane partitions) are 
combinatorial objects that are rather persistant against combinatorial reasonings.

\medskip

Our own proof of the alternating sign matrix theorem \cite{newproof} makes us believe that it could be helpful to work with signed enumerations: let $\alpha(n;k_1,\ldots,k_n)$ denote the number of monotone triangles with bottom row $k_1,\ldots,k_n$. The key identity in this 
proof is the following.
\begin{equation}
\label{identity}
\alpha(n;k_1,\ldots,k_n) = (-1)^{n-1} \alpha(n;k_2,\ldots,k_n,k_1-n)
\end{equation}
Obviously, this identity does not make any sense at first as $k_2, k_3, \ldots k_n, k_1-n$ is not strictly increasing if $k_1, k_2, \ldots, k_n$ is strictly increasing. However, it is not hard to see that, for fixed $n$,  the quantity
$\alpha(n;k_1,\ldots,k_n)$ can in fact be represented by a (unique) polynomial in $k_1,\ldots,k_n$ and so \eqref{identity}
can be understood as an identity for this polynomial. On the other hand, it is also possible to give $\alpha(n;k_1,\ldots,k_n)$ a combinatorial interpretation for all 
$(k_1,\ldots, k_n) \in \mathbb{Z}^n$ in terms of a signed enumeration. 
We have provided such an interpretation in \cite{simplified} and give three additional but related extensions in the concluding section of this article. These 
extensions provide combinatorial interpretations of \eqref{identity} and to give also a combinatorial proof of this identity could be an important step towards a combinatorial understanding of the alternating sign matrix theorem as we explain in Section~\ref{concluding}. It is hoped that a combinatorial proof of this identity as well as of other interesting identities involving monotone triangles follows the same lines as the combinatorial reasonings  we present in this article for Gelfand-Tsetlin patterns. 

\section{Definition of Gelfand-Tsetlin tree sequences}
\label{definition}

\begin{figure}
\begin{center}
\mbox{\scalebox{0.40}{%
    \includegraphics{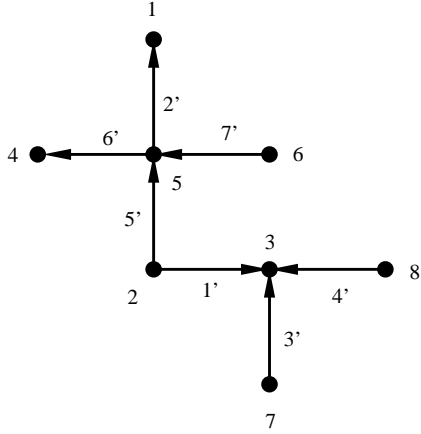}}}
\end{center}
\caption{\label{8tree} An $8$-tree.}   
\end{figure}

\begin{figure}
\begin{center}
\mbox{\scalebox{0.30}{%
    \includegraphics{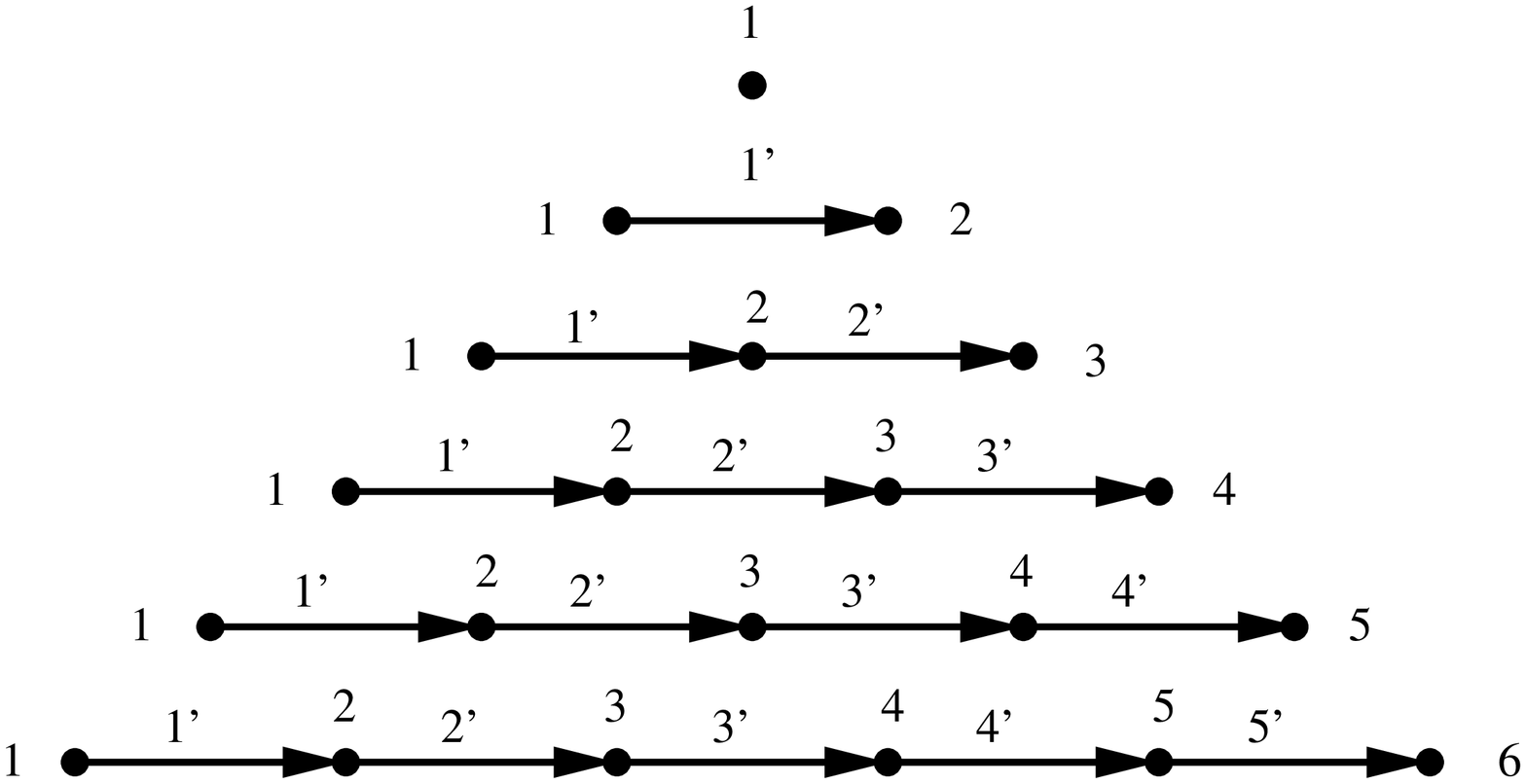}}}
\mbox{\scalebox{0.30}{%
    \includegraphics{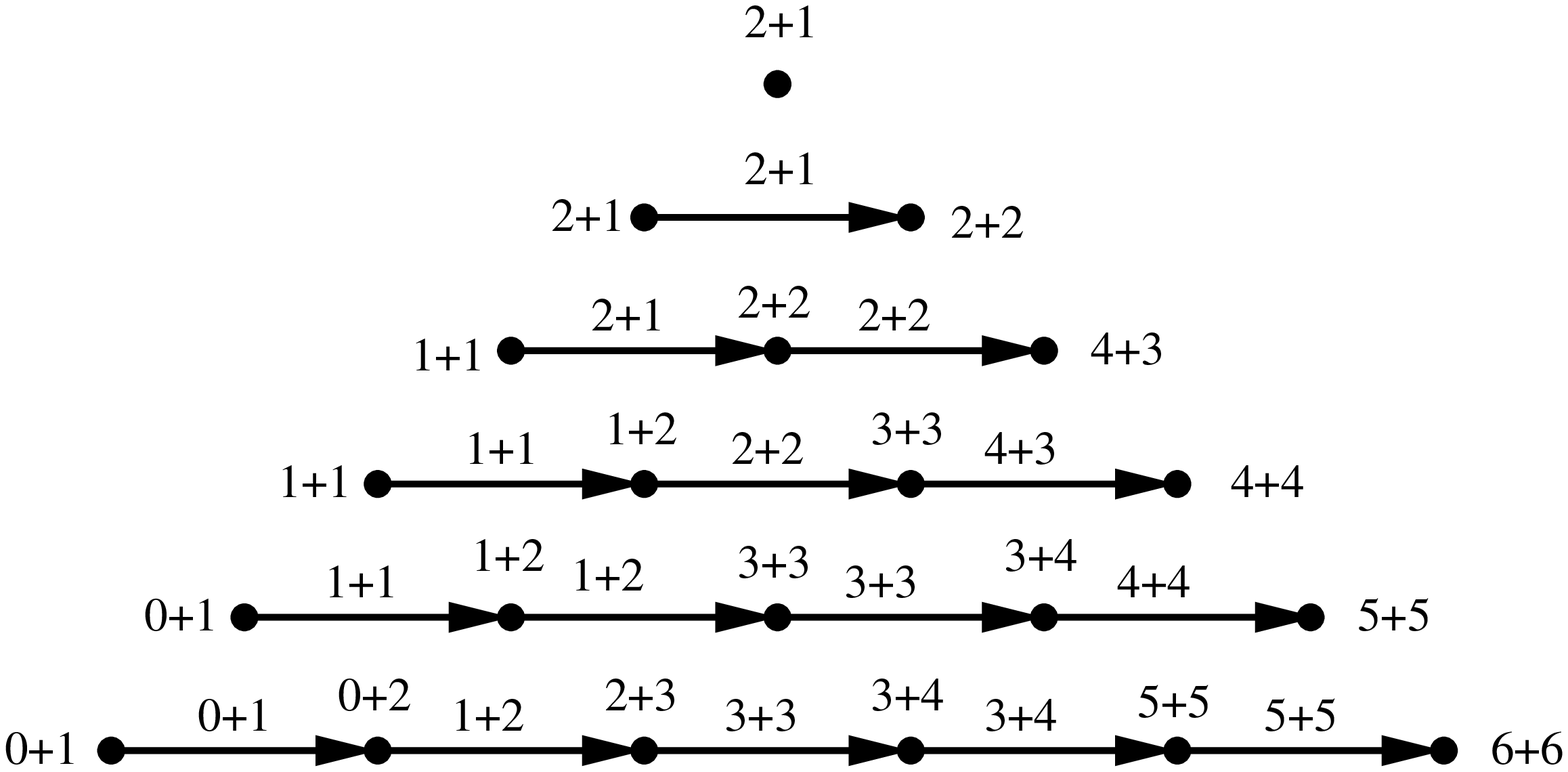}}}
\end{center}
\caption{Tree sequence for Gelfand-Tsetlin patterns of order $6$ and an example of an admissible labeling.}  
\label{gtpattern} 
\end{figure}

In
this paper, an {\bf $n$--tree} is a directed tree with $n$ vertices  such that the vertices are identified with integers in
$\{1,2,\ldots,n\}$ and the edges are identified with primed integers
in $\{1',2',\ldots,(n-1)'\}$. In Figure~\ref{8tree}, we give an example of an $8$--tree. We consider sequences of trees: a {\bf tree sequence of order 
$n$} is a sequence of trees ${\mathcal T}=(T_1,T_2,\ldots,T_n)$ 
such that $T_i$ is an $i$-tree for each $i$, see Figure~\ref{seq} for an example of order 
$5$. Each member of the family, the signed enumeration of which is given by \eqref{product}, will have a fixed underlying tree sequence of order $n$.  The actual objects will be certain admissible labelings (vertices and edges are labeled; the labels must not be confused with the ``names'' of the vertices and edges) of the underlying tree sequence. Gelfand-Tsetlin patterns will be one member of this family; in the underlying tree sequence ${\mathcal B} = (B_1,B_2,\ldots, B_n)$, 
the $i$-trees $B_i$ are paths with the canonial labeling, i.e. 
$j' = (j,j+1) \in E(B_i)$ for $j=1,2,\ldots, i-1$. In the following, the
tree $B_i$ will be referred to as the basic $i$-tree. In Figure~\ref{gtpattern}, we display the respective tree sequence of order 
$6$ (left figure) and the admissible labeling (a notion to be defined below) that corresponds to the Gelfand-Tsetlin pattern given in the introduction (right figure). In the right figure, we suppress the ``names'' of the vertices and edges in order to avoid a confusion with the labelings. However, these ``names'' are just the second summands of the labelings, whereas the first summand corresponds to the respective entry of the Gelfand-Tsetlin pattern given in the introduction.

\medskip

\begin{figure}
\begin{center}
\mbox{\scalebox{0.40}{%
    \includegraphics{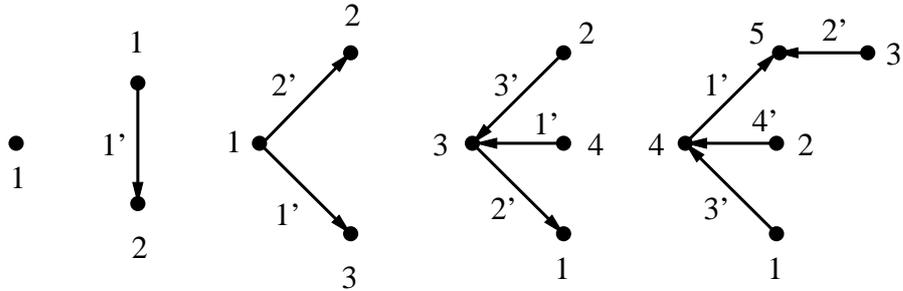}}}
\end{center}
\caption{\label{seq1} A tree sequence of order $5$.}   
\end{figure}

We work towards defining {\bf admissible labelings} of tree sequences.

\begin{defi}
Let $T$ be  an $n$--tree and ${\bf k}=(k_1,\ldots,k_n) \in \mathbb{Z}^n$.  
A vector ${\bf l}=(l_1,\ldots,l_{n-1}) \in \mathbb{Z}^{n-1}$ is said to be 
{admissible} for the pair $(T,{\bf k})$ if for each edge $j'=(p,q)$ of $T$ the following is fulfilled: if 
$k_p+p <  k_q+q$ then $k_p + p  \le l_j + j < k_q + q$ and otherwise $k_q+q \le  l_j+j < k_p+p$. In the latter case we say that the edge $j'$ is an inversion of the pair $(T,{\bf k})$.
\end{defi}

Phrased differently, if we 
label vertex $i$ with $k_i+i$ and edge $j'$ with $l_j+j$ for all $i$ and $j$ then, for each edge, the edge label is greater or equal than the minimum of the two vertex labels on the endpoints of the edge but smaller than the maximum. The edge is an inversion if it is directed from the maximum vertex label to the minimum vertex label. If, for an edge, the label of the tail coincides with the label of the head then there exists no vector $\bf l$ that is admissible for the pair $(T,{\bf k})$. In the following, we address the vectors ${\bf k}+(1,2,\ldots,n)$ and ${\bf l}+(1,2,\ldots,n-1)$ as the vertex labeling, respectively edge labeling of the tree and the vectors $\bf k$ and $\bf l$ as the shifted labelings.

\medskip

\begin{figure}
\begin{center}
\mbox{\scalebox{0.40}{%
    \includegraphics{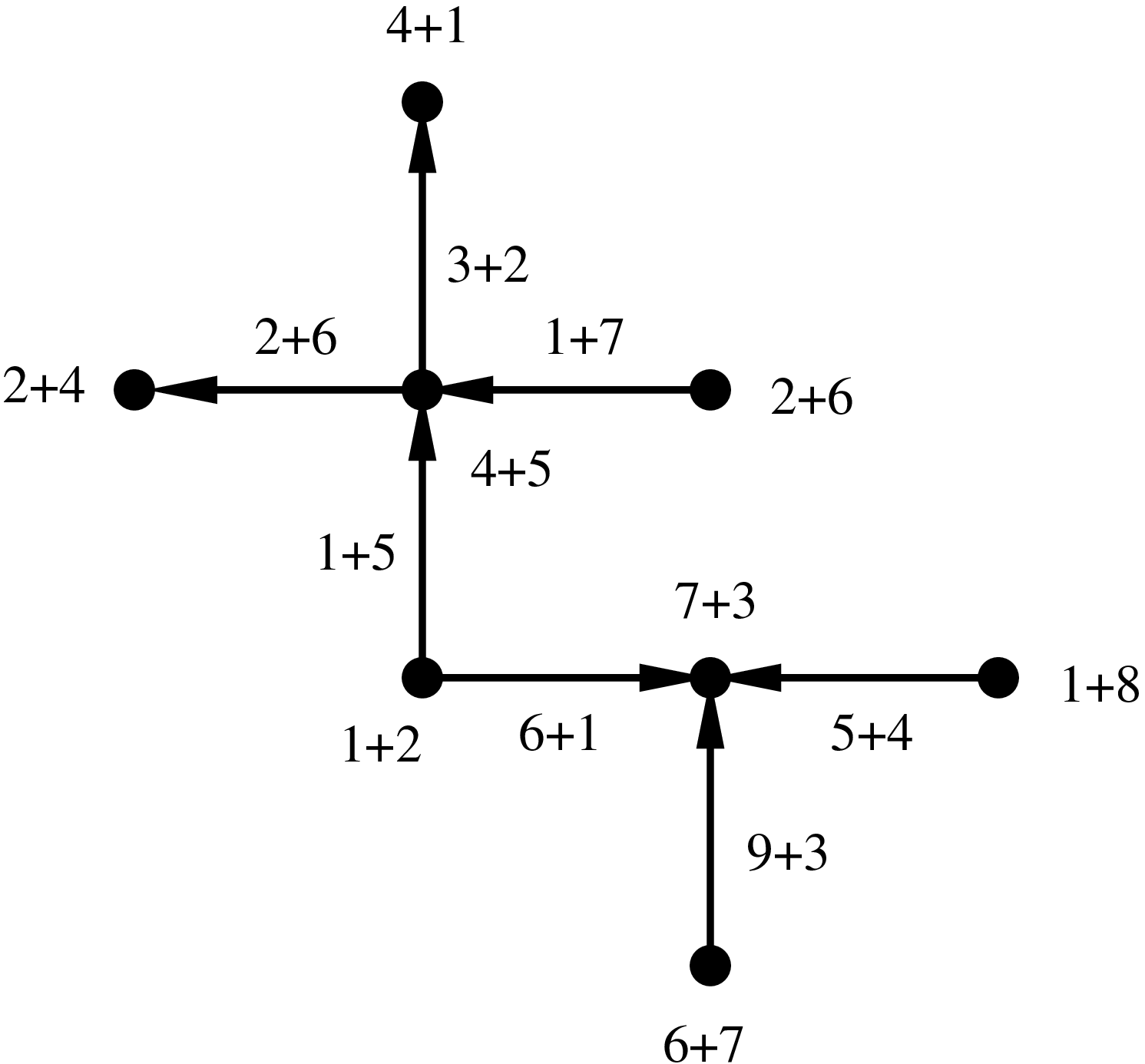}}}
\end{center}
\caption{\label{8treelabeling} An example of an admissible labeling.}   
\end{figure}

For instance, consider the $8$-tree $T$ in Figure~\ref{8tree} and the vector ${\bf k} = (4,1,7,2,4,2,6,1) \in \mathbb{Z}^8$. Then the vector
${\bf l}=(6,3,9,5,1,2,1)$ is admissible for $(T,{\bf k})$, see Figure~\ref{8treelabeling}. The
inversions are $2', 3', 6'$. Also observe that 
there is no admissible shifted labeling $\bf l$ if ${\bf k} = (4,1,7,2,4,2,6,2)$ as there is no $l_4$ with $2+8 \le l_4 + 4 < 7+3$.

\medskip

Now we are in the position to define  {\bf Gelfand-Tsetlin tree sequence}.

\begin{defi} A Gelfand-Tsetlin tree sequence
 associated with a tree sequence ${\mathcal T}=(T_1,\ldots,T_n)$  of order $n$ and a shifted labeling ${\bf k} \in \mathbb{Z}^n$ of the vertices of $T_n$ is a sequence $({\bf l_1, l_2, \ldots, l_n})$ of vectors ${\bf l_i} \in \mathbb{Z}^i$ with ${\bf l_n}= {\bf k}$ such that $\bf l_{i-1}$ is 
admissible for the pair $(T_i, {\bf l_i})$ if $i=2,3,\ldots,n$. We let ${\mathcal L}_n({\mathcal T},{\bf k})$ denote the set of these Gelfand-Tsetlin tree sequences. 
\end{defi}

\medskip

\begin{figure}
\begin{center}
\mbox{\scalebox{0.40}{%
    \includegraphics{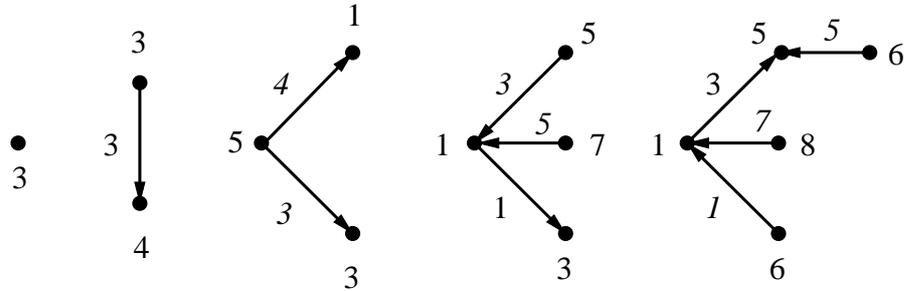}}}
\end{center}
\caption{A Gelfand-Tsetlin tree sequence.}   
\label{gttreesequence}
\end{figure}

In Figure~\ref{gttreesequence}, we give an example of a Gelfand-Tsetlin tree sequence 
associated with the tree sequence displayed in Figure~\ref{seq1}. Observe that 
${\bf k}=(5,6,3,-3,0)$ in this case. An edge label is displayed in italic type if the corresponding edge is an inversion. In Figure~\ref{gtpattern}, we represent the Gelfand-Tsetlin pattern from Section~\ref{intro}
as a Gelfand-Tsetlin tree sequence associated with $(B_1,B_2,B_3,B_4,B_5,B_6)$.

\medskip

We give a preliminary definition of the sign of a Gelfand-Tsetlin 
tree sequence:
the inversions of a Gelfand-Tsetlin tree sequence are the inversions
of the pairs $(T_i,{\bf l_i})$ for $i=2,3,\ldots,n$ and the sign is defined as
$(-1)^{\text{$\#$ of inversions}}$.
The (preliminary) sign of the Gelfand-Tsetlin tree sequence given in Figure~\ref{gttreesequence} is $-1$ as there are $7$ inversions. We will see that the signed enumeration of Gelfand-Tsetlin tree sequences associated with 
a fixed tree sequence ${\mathcal T}=(T_1,\ldots,T_n)$ of order $n$  and a fixed shifted labeling ${\bf k}=(k_1,\ldots,k_n)$ of the vertices of $T_n$ is, up to a sign, equal to \eqref{product}.
This sign only depends on the underlying unlabeled tree sequence ${\mathcal T}$ and will be defined next. After that we adjust the definition of the sign of a Gelfand-Tstelin tree sequence by multiplying this global sign.

\medskip

For this purpose, we define the {\bf sign of an $n$--tree} $T$: fix a
{\it root vertex} $r$ of the tree. The {\it standard orientation} with
respect to this root is the orientation in which each edge is
oriented away from the root. An edge in $T$ is said to be a {\it
reversed} edge if its orientation does not coincide with the
standard orientation.  If, in our example in Figure~\ref{8tree}, we choose $2$ to be the root then the reversed edges are $3'$, $4'$ and $7'$. Except for the root, each vertex is the head
of a unique edge with respect to the standard orientation. We
obtain a permutation $\pi$ of $\{1,2,\ldots,n\}$, if we order the
head vertices of the edges in accordance with their edge names (i.e.
for the edges  $i'=(a,b)$ and $j'=(c,d)$ with $i < j$, the vertex $b$ comes before
vertex $d$ in the permutation) and prepend the root $r$ at the
beginning of the permutation. In our running example, we obtain the 
permutation $\pi= 2  \, 3 \, 1 \, 7 \, 8 \, 5 \, 4 \, 6$. Then the sign of $T$ is defined as
follows.
\begin{equation}
\label{sign} \sgn T = (-1)^{\text{$\#$ of reversed edges}} \sgn \pi
\end{equation}
The sign of the tree in Figure~\ref{8tree} is $1$ as there are $3$ reversed edges and 
$\sgn \pi = -1$.

\medskip

We need to show that the sign does not depend on the choice of the root: 
suppose $s$ is
a vertex adjacent to the root $r$. If we change from root $r$ to
root $s$, we have to interchange $r$ and $s$ in the permutation
$\pi$, which reverses the sign of $\pi$. This is because the
standard orientation with respect to the root $s$ coincides with
the standard orientation with respect to the root $r$ except for
the edge incident with $r$ and $s$, where the orientation is
reversed. For the same reason, shifting the root from $r$ to $s$, either
increases or decreases the number of reversed edges by $1$. 
Consequently, the product in \eqref{sign} remains unaffected.

\medskip

The {\bf sign of a tree sequence} ${\mathcal T}=(T_1,T_2,\ldots,T_n)$ is defined as the product of the signs of 
the $i$-trees in the sequence, i.e.
$$
\sgn {\mathcal T} = \sgn T_1 \cdot \sgn T_2 \cdots \sgn T_n.
$$
The sign of the tree sequence in Figure~\ref{seq1} is $-1$ as $\sgn T_1=1, 
\sgn T_2=1,  \sgn T_3=-1, \sgn T_4=1,$ and $\sgn T_5=1$. Concerning Gelfand-Tsetlin patterns we 
obviously have $\sgn B_i = 1$, which implies $\sgn {\mathcal B} = 1$. 

\medskip

Here is the final definition of the {\bf sign  of a Gelfand-Tsetlin tree sequence} 
${\bf L}=({\bf l_1, l_2, \ldots, l_n}) \in {\mathcal L}_n({\mathcal T},{\bf k})$: 
$$
\sgn {\bf L}  =  (-1)^{\text{$\#$ of inversions of 
${\bf L }$}} \cdot
\sgn {\mathcal T}
$$
The signed enumeration of elements in $ {\mathcal L}_n({\mathcal T},{\bf k})$ is denote by $L_n({\mathcal T}, {\bf k})$. The sign of the Gelfand-Tsetlin tree sequence given in Figure~\ref{gttreesequence} is $1$ as there are $7$ inversions and the sign of the underlying unlabeled tree sequence is $-1$. We are in the position to state an important result of this paper.

\medskip

\begin{theo} 
\label{main}
The signed enumeration  of
Gelfand-Tsetlin tree sequences associated with a fixed underlying unlabeled tree sequence 
${\mathcal T}=(T_1,\ldots,T_n)$ of order $n$  and a shifted labeling ${\bf k}=(k_1,\ldots,k_n)$ of the vertices of $T_n$ is  given by
$$
\prod_{1 \le i < j \le n} \frac{k_j-k_i+j-i}{j-i}.
$$
\end{theo}

\medskip

Before we turn our attention to searching for properties of $L_n({\mathcal T}, {\bf k})$ that determine this quantity uniquely, we want to mention an obvious generalization of Gelfand-Tsetlin tree sequences, which  we do not consider in this article, but might be interesting to look at: the notion of admissibility makes perfect sense if the tree $T$ is replaced by any other graph. Are there any nice assertions to be made on ``Gelfand-Tsetlin graph sequences''?

\section{Properties of $L_n({\mathcal T}, {\bf k})$: independency and shift-antisymmetry}
\label{indsym}

We say that a function $f(k_1,\ldots,k_n)$ on $\mathbb{Z}^n$ is  {\bf shift-antisymmetric} iff 
$$
f(k_1,\ldots,k_n)=-f(k_1,\ldots,k_{i-1},k_j+j-i,k_{i+1},\ldots,k_{j-1},k_i+i-j,k_{j+1},\ldots,k_n)
$$
for all $i,j$ with $1 \le i < j \le n$ and all $(k_1,\ldots,k_n) \in \mathbb{Z}^n$.  In this section we prove by induction with respect to $n$ that 
the signed 
enumeration $L_n({\mathcal T}, {\bf k})$  has the following two properties. 

\medskip

\begin{itemize}
\item {\bf Independency:}  $L_n({\mathcal T}, {\bf k})$  
does not depend on the tree sequence ${\mathcal T}$.  
\item {\bf Shift-antisymmetry:}  
$L_n({\mathcal T}, {\bf k})$ is 
shift-antisymmetric in ${\bf k}=(k_1,\ldots,k_n)$. In fact, we prove the following stronger result: fix $i,j$ with $1 \le i < j \le n$. We construct a tree sequence of order $n$, denoted by 
${\mathcal S}_{n}^{i,j}$, and an associated sign reversing involution on the set of Gelfand-Tsetlin tree sequences of the tree sequence ${\mathcal S}^{i,j}_n$ such that the shifted vertex labeling ${\bf k} \in \mathbb{Z}^n$ of 
the largest tree is transformed into 
$$ E^{j-i}_{k_{j}} E^{i-j}_{k_i} S_{k_i,k_{j}} {\bf k} = 
(k_1,\ldots,k_{i-1},k_j+j-i,k_{i+1},\ldots,k_{j-1},k_{i}+i-j,k_{j+1},\ldots,k_n),$$
where $S_{x,y} f(x,y) = f(y,x)$ and $E_x p(x) = p(x+1)$.
\end{itemize}

\medskip

The proofs are combinatorially in the following sense: suppose we are given two sets $A$ and $B$ and a signed enumeration $|.|_{-}$ on each of the sets such that 
$|A|_{-}=|B|_{-}$. Then we find decompositions of $A$ and $B$ into two 
sets $A_1, A_2$ and $B_1,B_2$, respectively,  such that there is a 
sign preserving bijection between 
$A_1$ and $B_1$ and $|A_2|_{-}=|B_2|_{-}=0$, where the latter identities are proven by giving 
sign reversing involutions on $A_2$ and $B_2$. However, if we have  $|A|_{-}=-|B|_{-}$ then 
the bijection between $A_1$ and $B_1$ is sign reversing.

\medskip

Observe that there is nothing to prove for $n=1$. 
We deal with the independency first.
\begin{lem}
\label{invariant}
The independency and shift-antisymmetry for order $n-1$ implies the independency for order $n$.  
\end{lem}

\medskip

{\it Proof.} For a tree sequence ${\mathcal T}=(T_1,T_2,\ldots,T_n)$ of order 
$n$ we have  
$$
L_n({\mathcal T},{\bf k}) = \sgn T_n \cdot 
(-1)^\text{$\#$ of inversions of $(T_n,{\bf k})$}
\sum_{\text{${\bf l} \in \mathbb{Z}^{n-1}$ is admissible for
$(T_n,{\bf k})$}} 
L_{n-1}({\mathcal T}_{<n}, \bf l),
$$
where ${\mathcal T}_{<n} = (T_1,T_2,\ldots,T_{n-1})$.
The independency for $n-1$ implies that $L_n({\mathcal T},{\bf k})$
is invariant under the replacement of ${\mathcal T}_{<n}$ by any other 
tree sequence of order $n-1$. We have to show that it is also 
invariant under the replacement of $T_n$ by any other $n$-tree.
The strategy  is as follows: we first show that 
$L_n({\mathcal T},{\bf k})$ is invariant under certain tree 
operations on $T_n$ and then verify that every tree can be obtained from every other by means of these operations. To prove this invariance, we often
replace ${\mathcal T}_{<n}$ by a particularly convenient tree sequence.

\medskip

We define the first tree operation: let $T_n$ be an $n$-tree and $T'_n$ be an $n$-tree which is obtained from $T_n$ by reversing the orientation of a single edge. Then 
$\sgn T_n =  - \sgn T'_n$,  the number of inversions of $(T_n, {\bf k})$ differs from the number of inversions of $(T'_n, {\bf k})$ by $1$ and 
${\bf l} \in \mathbb{Z}^{n-1} $ is admissible for $(T_n, {\bf k})$ if and only if ${\bf l}$ is admissible for $(T'_n, {\bf k})$. This implies that 
$L_n({\mathcal T},{\bf k})$ is invariant under the replacement of 
$T_n$ by $T'_n$. 

\medskip

\begin{figure}
\begin{center}
\mbox{\scalebox{0.40}{%
    \includegraphics{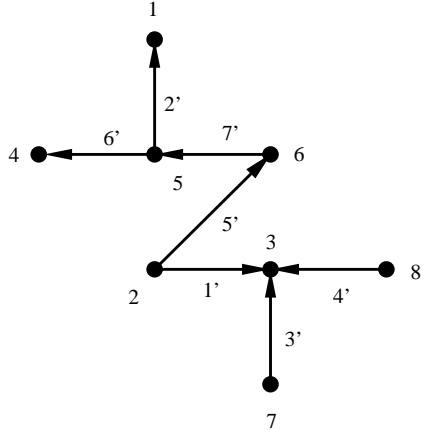}}}
\end{center}
\caption{\label{sliding} The $8$-tree is obtained from the $8$-tree in Figure~\ref{8tree} by sliding edge $5'$ along edge $7'$. }   
\end{figure}

\medskip

For the second operation we assume $n \ge 3$. It is illustrated in Figure~\ref{sliding} and defined as follows: suppose that $i'$ and $j'$ are two edges in the $n$-tree $T_n$ that have a vertex $q$ in common. Let $T'_n$ be the tree we obtain from $T_n$ by replacing vertex $q$ in $i'$ by the vertex of $j'$ which is different from $q$.  Then we say that $T'_n$ is obtained from $T_n$ by {\bf sliding edge $i'$ along edge $j'$}. In the following argument, we let $p$ be the vertex of $i'$ in $T_n$ that is different from $q$ and $r$ be the vertex of $j'$ that is different from $q$.

\medskip

We show $\sgn T_n = \sgn T'_n$:  let $q$ be the root. The head of the old edge $i'$ (i.e. in $T_n$) as well as of the new edge $i'$ (i.e. in $T'_n$) is $p$ with respect to the standard orientation. Moreover, the edge $i'$ is reversed in $T_n$ if and only if it is reversed in $T'_n$. 
There is no change for the remaining edges, since the standard orientation 
does not change there. Hence, neither the permutation $\pi$ nor the set of 
reversed edges is changed.

\medskip

In order to show that $L_n({\mathcal T},{\bf k})$ is invariant under the replacement of $T_n$ by $T'_n$, we have to distinguish between the six possibilities for the relative positions of 
$k_p+p, k_q+q, k_r+r$. As we have a symmetry between vertex $q$ and vertex $r$  we may assume without loss of generality that $k_q+q \le k_r +r$. We let 
${\mathcal T}'$ denote the tree sequence that we obtain from $\mathcal T$
by replacing $T_n$ by $T'_n$.

\medskip 

{\it Case 1.} $k_p + p \le k_q +q \le k_r +r$: we decompose 
${\mathcal L}_n({\mathcal T}',{\bf k})$ into two sets as follows. Let 
${\bf l} \in \mathbb{Z}^{n-1}$ be an admissible shifted edge labeling
of $T'_n$. The first set contains the Gelfand-Tsetlin tree sequences where 
the label of edge $i'$ fulfills 
$k_p + p \le l_i + i < k_q + q$, whereas for the second set we have  
$k_q + q \le l_i + i \le k_r+r$. The signed 
enumeration of the first set is obviously equal to 
${\mathcal L}_n({\mathcal T},{\bf k})$, since the edge $i'$ is an inversion
of
$T_n$ if and only if it is an inversion of $T'_n$. We have to show that the signed enumeration 
of the second set reduces to zero: we replace ${\mathcal T}_{< n}$ by 
${\mathcal S}^{i,j}_{n-1}$. As $k_q + q \le l_i + i < k_r + r$ and 
$k_q + q \le l_j + j < k_r +r$, the sign reversing involution on 
the set of all Gelfand-Tsetlin tree sequence associated with 
${\mathcal S}^{i,j}_{n-1}$ induces a sign reversing involution 
on the second subset of ${\mathcal L}_n({\mathcal T}',{\bf k})$.

\medskip

{\it Case 2.} $k_q +q \le k_p + p \le k_r +r$: if
${\bf l} \in \mathbb{Z}^{n-1}$ is an admissible shifted edge labeling
of $T_n$ for an element of ${\mathcal L}_n({\mathcal T},{\bf k})$ then we 
have $k_q +q \le l_i + i < k_p + p$; in  ${\mathcal L}_n({\mathcal T}',{\bf k})$ we have $k_p + p \le l_i + i < k_r +r$. The edge $i'$ is an inversion for the pair $(T_n,{\bf k})$ if and only it is no inversion for the pair 
$(T'_n,{\bf k})$. We decompose both sets into two sets according to the edge label of $j'$: in the 
first set we have 
$k_q+q \le l_j + j < k_p + p$ and in the second set we have $k_p + p \le l_j + j < k_r +r$. If we replace ${\mathcal T}_{< n}$ by ${\mathcal S}_{n-1}^{i,j}$, we see that in case of ${\mathcal L}_n({\mathcal T},{\bf k})$ the signed enumeration of the first set is zero, while for ${\mathcal L}_n({\mathcal T}',{\bf k})$ the signed enumeration of the second set is zero. For the two other sets, the replacement of $(l_i,l_j) \to (l_j+j-i,l_i+i-j)$ of the shifted edge labels of the largest tree and performing the sign reversing involution on ${\mathcal S}_{n-1}^{i,j}$ is a sign preserving involution.

\medskip

{\it Case 3.} $k_q + q \le k_r +r \le k_p + p$: for the edge label of $i'$ in 
$T_n$ we have $k_q + q \le l_i + i < k_p + p$. We decompose ${\mathcal L}_n({\mathcal T},{\bf k})$ into two sets, where we have $k_{q}+q \le l_i +i < k_r + r$ and 
$k_r + r \le l_i + i < k_p + p$, respectively. As $k_q + q \le l_j + j < k_r + r$, the signed enumeration of the first set is zero, while the signed enumeration of the 
second set coincides with the signed enumeration of the elements in 
${\mathcal L}_n({\mathcal T}',{\bf k})$.  

\medskip 

In order to conclude the proof of Lemma~\ref{invariant}, it suffices to show that every $n$-tree can be transformed 
into every other by means of the two operations ``sliding an edge along another edge'' and ``reversing the orientation of an edge''.  As both operations are in fact involutions, it suffices to show that every $n$-tree can be transformed into the basic $n$-tree $B_n$. First of all, it is obvious that sliding and reversing can be used to transform a given $n$-tree into a directed path. Hence, it suffices to show that it is possible to interchange 
vertices as well as edges. In both cases, it suffices to consider 
adjacent vertices, respectively
edges. Concerning edges, suppose $x'$ and $y'$ are adjacent edges. By possibly reversing the 
orientation of one edge, we may assume without loss of generality that $x'=(a,b)$ and $y'=(b,c)$. Then 
the following sequence of operations interchanges  the edges:
\begin{multline*}
x'=(a,b), y'=(b,c) \rightarrow x'=(a,c), y'=(b,c) \rightarrow 
x'=(a,c), y'=(b,a) \\ \rightarrow x'=(b,c), y'=(b,a) \rightarrow 
x'=(b,c), y'=(a,b)
\end{multline*}
(Note that all operations except for the last are slides, which implies that 
interchanging edges reverses the sign of the $n$-tree.)
Concerning swapping vertices, assume that we want to interchange vertex $a$ and $b$ and that $x'=(a,b)$ is an 
edge. We reverse the orientation of $x'$  and slide all edges incident with $a$ but different from $x'$ along $x'$ to $b$ as well as all edges incident with $b$ but different from $x'$ along $x'$ to $a$. (Again we see that swapping vertices reverses the sign.)
 \qed

\medskip

Now we turn to the shift-antisymmetry.

\begin{lem}
\label{antisym}
 The independency for order $n$ implies the shift-antisymmetry for 
order $n$.
\end{lem}

\medskip

\begin{figure}
\begin{center}
\mbox{\scalebox{0.40}{%
    \includegraphics{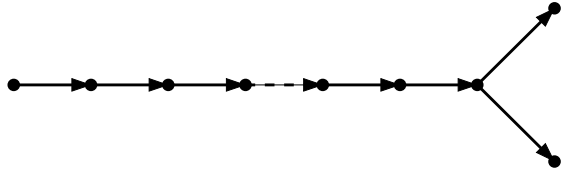}}}
\end{center}
\caption{\label{sym} The tree $S_m$.}   
\end{figure}

{\it Proof.} Fix $i,j$ with $1 \le i < j \le n$. We define a tree sequence 
${\mathcal S}^{i,j}_n=(T_1,\ldots,T_n)$ of order $n$: let $S_m$ be the directed tree with $m$ vertices sketched in Figure~\ref{sym} and, for $3 \le m \le n$,
let this be the underlying tree for $T_m$. (Note that there is no choice for the underlying tree if $m=1,2$.) There are no restrictions on the names of the vertices and edges except that the two sinks in $T_{n}$ are $i$ and $j$, the two sinks in $T_{n-1}$ are the unprimed versions of the 
edges incident with $i$ and $j$ in $T_{n-1}$, the two sinks in $T_{n-2}$ are the unprimed versions of the edges incident with the two sinks in $T_{n-1}$ etc. Let ${\bf k}=(k_1,\ldots,k_{n}) \in \mathbb{Z}^n$ and ${\bf k'} = E^{j-i}_{k_j} E^{i-j}_{k_i}S_{k_i,k_j} {\bf k}$. Then the following is a sign reversing involution between the Gelfand-Tsetlin tree sequence associated with ${\mathcal S}^{i,j}_n$ and fixed shifted vertex labeling $\bf k$ of $T_n$ and those where the shifted vertex labeling of  $T_n$ is given by $\bf k'$: for $m \ge 3$, we interchange in $T_m$ the labels of the two sink vertices as well as the labels of the two edges incident with the sinks; in $T_2$ we interchange the two vertex labels. This either produces or resolves an inversion in $T_2$ and concludes the proof of Lemma~\ref{antisym}. 

\medskip 

Alternatively, we can also argue as follows: let $T'_n$ be the tree which we obtain from $T_n$ by interchanging vertex $i$ and vertex $j$ (the underlying tree remains unaffected) and ${\mathcal T}'=(T_1,\ldots,T_{n-1},T'_{n})$.  As $\sgn T_n = - \sgn T'_n$, we obviously have 
$$
L_n({\mathcal T},{\bf k}) = - E^{j-i}_{k_j} E^{i-j}_{k_i}S_{k_i,k_j} L_n({\mathcal T}',{\bf k}).
$$
The assertion follows from Lemma~\ref{invariant} since $L_n({\mathcal T}',{\bf k})=L_n({\mathcal T},{\bf k})$. \qed

\medskip 

In Appendix~\ref{shiftantisymmetry}, a direct combinatorial proof of the shift-antisymmetry of the enumeration formula for Gelfand-Tsetlin patterns is sketched, which does not make use of the notion of Gelfand-Tsetlin tree sequences.

\section{Taking differences --  
${L}_n(\mathcal T, {\bf k})$ is a polynomial}
\label{diff}

The quantity ${L}_n(\mathcal T, {\bf k})$ is not characterized by the properties we have derived so far.  Next, we show that ${L}_n(\mathcal T, {\bf k})$ is a polynomial of 
degree no greater than $n-1$ in every $k_i$, which is the last ingredient to finally see that 
it is equal to \eqref{product}.

\medskip

In order to show that $p(x)$ is a polynomial in $x$ of degree no greater than $n-1$,  it suffices to prove that $\Delta^n_x p(x) = 0$ where $\Delta_x:=E_x - \id$ is the difference operator.
Thus it suffices to show the following. 

\begin{lem}
\label{degreeequation}
For $i \in \{1,2,\ldots,n\}$ we have   $\Delta^n_{k_i} {L}_n({\mathcal T},{\bf k}) = 0$ .
\end{lem}

\medskip

\begin{figure}
\begin{center}
\mbox{\scalebox{0.40}{%
    \includegraphics{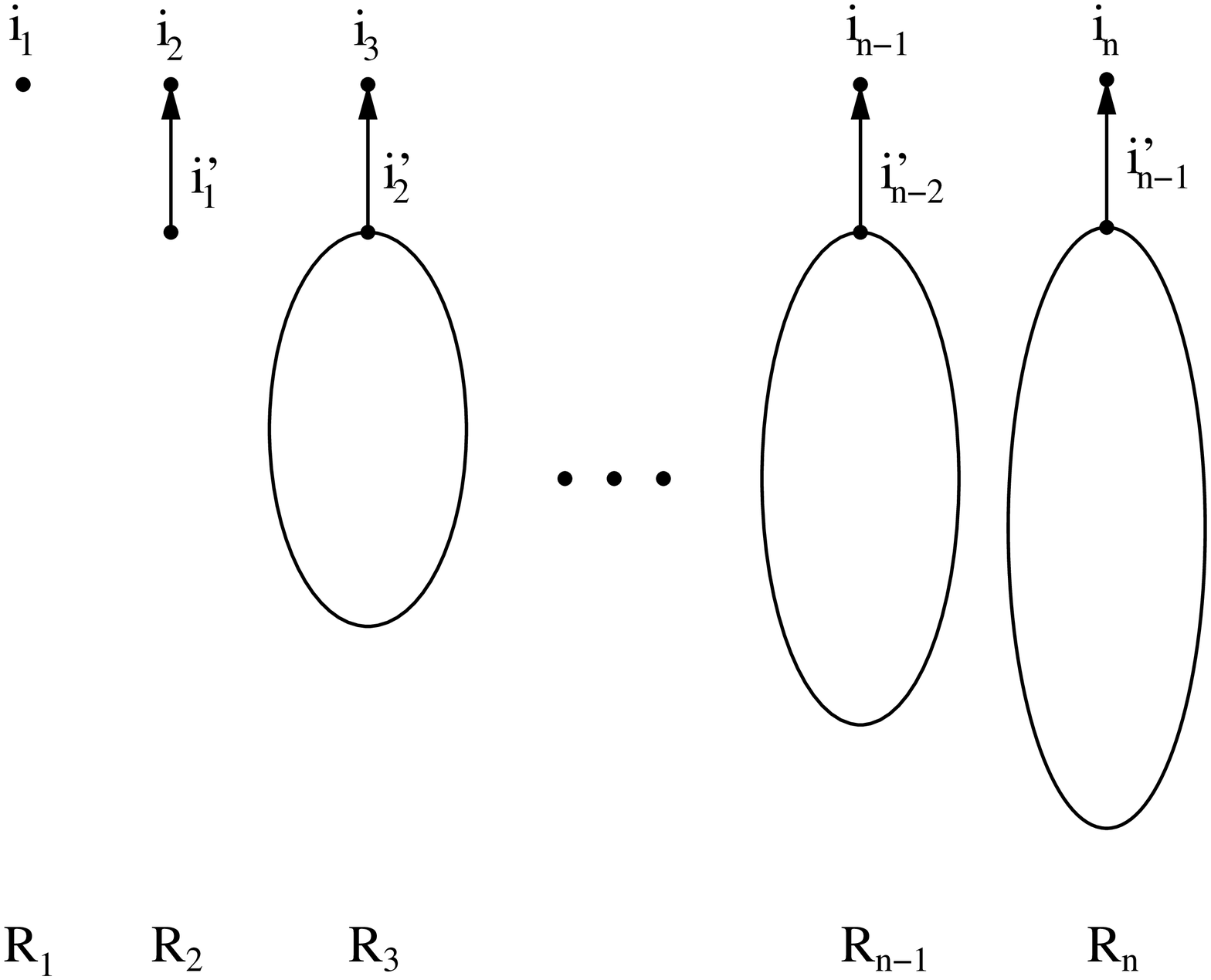}}}
\end{center}
\caption{\label{degreefigure} Tree sequence in the proof of Lemma~\ref{degreeequation}.}   
\end{figure}

{\it Proof.} 
We define a convenient tree sequence ${\mathcal R_{n,i}} = (R_1,\ldots, R_n)$ (see Figure~\ref{degreefigure}) and find a combinatorial interpretation for $\Delta^j_{k_i} {L}_n({\mathcal R_{n,i}},{\bf k})$ if $j \in \{0,1,\ldots,n-1\}$: in $R_n$, we require $i=:i_n$ to be a leaf, in $R_{n-1}$ we require the unprimed version $i_{n-1}$ of the edge incident with $i_n$ in $R_n$ to be a leaf, in $R_{n-2}$ we require the unprimed version $i_{n-2}$ of the edge incident with $i_{n-1}$ in 
$R_{n-1}$ to be a leaf etc. As for the orientations of the edges $i'_1, i'_2,\ldots,i'_{n-1}$, we choose the vertices $i_2,i_3,\ldots,i_n$ to be sinks. By $l_{i_1}+i_1,l_{i_2}+i_2,\ldots, l_{i_{n-1}}+i_{n-1}$, we denote the respective edge labels (which are of course also vertex labels in the next level).  

\medskip

We define $\Delta^j_{k_i} {\mathcal L}_n({\mathcal R_{n,i}},{\bf k})$: it is 
the set of labeled tree sequences on the unlabeled tree sequence ${\mathcal R}_{n,i}$ such that the conditions 
on the edge labels are as for Gelfand-Tsetlin tree sequence in 
${\mathcal L}_n({\mathcal R_{n,i}},{\bf k})$, except  
for the edges $i'_{n-j},i'_{n-j+1},\ldots,i'_{n-1}$ in $R_{n-j+1}, R_{n-j+2}, \ldots, R_n$, respectively, where we require
$l_{i_{n-j}}+i_{n-j}=l_{i_{n-j+1}}+i_{n-j+1}=\ldots=l_{i_{n-1}}+i_{n-1}=k_i+i$. 
As for the sign, we compute it as usual only we ignore the contributions
of the edges $i'_{n-j} \in E(R_{n-j+1}),i'_{n-j+1} \in E(R_{n-j+2}),\ldots,i'_{n-1} \in E(R_n)$.

\medskip

Then, by induction with respect to $j$, the signed enumeration of these labeled tree sequences on ${\mathcal R}_{n,i}$ is equal to $\Delta^j_{k_i} {L}_n({\mathcal R_{n,i}},{\bf k})$: for $j=0$ this is obvious. It suffices to show that
$$
\Delta_{k_i} | \Delta_{k_i}^j {\mathcal L}_n({\mathcal R}_{n,i},{\bf k})|_{-} 
= |\Delta_{k_i}^{j+1} {\mathcal L}_n({\mathcal R}_{n,i},{\bf k})|_{-}.
$$ Consider an element from $E_{k_i} \Delta^j_{k_i} {\mathcal L}_n({\mathcal R_{n,i}},{\bf k})$ such that the vertex label  of the sink $i_{n-j}$ of the edge $i'_{n-j-1}$ in $R_{n-j}$ (which is $l_{i_{n-j}}+i_{n-j}=k_i+i+1$) is greater than the vertex label of the other endpoint of the edge. Then, by decreasing the labels $l_{i_{n-j}}+i_{n-j},l_{i_{n-j+1}}+i_{n-j+1},\ldots, l_{i_{n-1}}+i_{n-1}, k_i+i+1$ (which are all equal) by $1$, we obtain a corresponding element in $\Delta^j_{k_i} {\mathcal L}_n({\mathcal R_{n,i}},{\bf k})$, except for the case when 
$l_{i_{n-j-1}}+i_{n-j-1}=k_{i}+i$. In such a tree sequence, we also decrease the labels $l_{i_{n-j}}+i_{n-j},l_{i_{n-j+1}}+i_{n-j+1},\ldots, l_{i_{n-1}}+i_{n-1}, k_i+i+1$ by $1$ to obtain an element of  $\Delta^{j+1}_{k_i} {\mathcal L}_n({\mathcal R_{n,i}},{\bf k})$. This way, we obtain exactly the elements of $\Delta^{j+1}_{k_i} {L}_n({\mathcal R_{n,i}},{\bf k})$ such that the edge $i'_{n-j-1}$ is no inversion in $R_{n-j}$.
On the other hand, if the edge $i'_{n-j-1}$ is an inversion for an element of
$\Delta^j_{k_i} {\mathcal L}_n({\mathcal R_{n,i}},{\bf k})$, then, by increasing the labels $l_{i_{n-j}}+i_{n-j},l_{i_{n-j+1}}+i_{n-j+1},\ldots, l_{i_{n-1}}+i_{n-1}, k_i+i$ by $1$, we obtain a corresponding element in 
$E_{k_i} \Delta^j_{k_i} {\mathcal L}_n({\mathcal R_{n,i}},{\bf k})$, except for the case when $l_{i_{n-j-1}}+i_{n-j-1}=k_{i}+i$. This way, we obtain exactly the elements of $\Delta^{j+1}_{k_i} {L}_n({\mathcal R_{n,i}},{\bf k})$ such that the edge $i'_{n-j-1}$ is an inversion in $R_{n-j}$. The sign that comes from the inversion $i'_{n-j-1}$ in $R_{n,n-j}$ takes into account for the fact that we ``subtract'' the greater set from the smaller set in this case. 

\medskip

Now observe that in fact $\Delta^{n-1}_{k_i} {\mathcal L}_n({\mathcal R_{n,i}},{\bf k})$ does not depend on $k_i$ and, consequently,  $\Delta^{n}_{k_i} {L}_n({\mathcal R_{n,i}},{\bf k})$ must be zero. \qed

\medskip

\medskip

We are finally in the position of prove Theorem~\ref{main}.

\medskip

{\it Proof of Theorem~\ref{main}.}  By the shift-antisymmetry (Lemma~\ref{antisym}), we conclude that the polynomial (Lemma~\ref{degreeequation}) ${L}_n({\mathcal T},{\bf k})$ vanishes if $k_i+i=k_j+j$ for distinct $i,j \in \{1,2,\ldots,n\}$. This implies that 
the expression in \eqref{product} has to be a factor of ${L}_n({\mathcal T},{\bf k})$. Again by Lemma~\ref{degreeequation}, we know that it is a polynomial of degree no geater than $n-1$ and since \eqref{product} is of degree $n-1$ in every $k_i$, this implies that 
$$
{L}_n({\mathcal T},{\bf k}) = C \cdot \prod\limits_{1 \le i < j \le n} \frac{k_j - k_i + j - i}{j-i}, 
$$
where $C \in \mathbb{Q}$.
As there is only one Gelfand-Tsetlin pattern with bottom row $(1,1,\ldots, 1) \in \mathbb{Z}^n$, we can conclude that $C=1$.

\medskip

The combinatorial interpretation of $\Delta^j_{k_i} {\mathcal L}_n({\mathcal R_{n,i}},{\bf k})$ was surely the main ingredient in the proof of Lemma~\ref{degreeequation}. The remainder of this section is devoted to use basically the same idea to give a combinatorial proof of the identity
\begin{equation}
\label{toshow}
e_{\rho}(\Delta_{k_1},\ldots,\Delta_{k_n}) {L}_n({\mathcal T},{\bf k}) = 0,
\end{equation}
which holds for ${\rho} \ge 1$ and where 
$$
e_{\rho}(X_1,\ldots,X_n) = \sum_{1 \le i_1 < i_2 < \ldots < i_{\rho} \le n} X_{i_1} X_{i_2} \cdots X_{i_{\rho}}
$$
is the $\rho$-th elementary symmetric function. (An algebraic proof, which already uses the fact that 
$L_n({\mathcal T},{\bf k})$ is equal to \eqref{product} as well as the presentation of 
\eqref{product} in terms of a determinant (see \eqref{det}), can be found in  \cite[Lemma~1]{newproof}.) This identity is of interest as it is the crucial fact in the proof of \eqref{identity} given 
in \cite{newproof}. 

\medskip

Even though the ideas are straight forward, this combinatorial proof of \eqref{toshow} is a bit elaborate. (However, nothing else is to be expected when a statement is related to alternating sign matrix counting.) In fact, the benefit of 
this exercise is not primarily the proof of \eqref{toshow} but an improvement of the understanding of how to 
interpret the application of difference operators to enumerative quantities such as 
${L}_{n}({\mathcal T},{\bf k})$ combinatorially. To give a hint as to why such an understanding could be of interest, observe that the proof of \eqref{toshow} relies on a combinatorial interpretation of 
\begin{equation}
\label{subset}
\Delta_{k_{i_1}} \Delta_{k_{i_2}} \ldots \Delta_{k_{i_\rho}} 
{L}_{n}({\mathcal T},{\bf k})
\end{equation}
for subsets $\{i_1,\ldots,i_{\rho}\} \subseteq [n]$. As the number of monotone triangles with bottom row
$(k_1,\ldots,k_n)$ is given by 
\begin{equation}
\label{operator2}
\alpha(n;k_1,\ldots,k_n) = \left( \prod_{1 \le p < q \le n} (\id + 
\Delta_{k_p} \delta_{k_q}) \right)  {L}_{n}({\mathcal T},{\bf k}), 
\end{equation}
where $\delta_x = \id - E^{-1}_x$ is a second type of difference operator (see 
Section~\ref{concluding}), ideas along these lines might also lead to a combinatorial proof of this formula. 

\medskip

We need a more general notion of 
admissibility. The idea is simple and very roughly as follows: we require each vertex of a fixed vertex set $R$ of the tree $T$ to have an associated edge incident with it such that the edge label takes on the extreme label given by the vertex label. 

\medskip

\begin{defi}
\label{weakad}
Given an $n$-tree $T$, an $n$-tupel ${\bf k}=(k_1,\ldots,k_n) \in \mathbb{Z}^{n}$ and a subset $R \subseteq [n]=:\{1,2,\ldots,n\}$  of vertices of $T$,  we define a vector ${\bf l}=(l_1,\ldots,l_{n-1}) \in 
\mathbb{Z}^{n-1}$  to be 
{\bf weakly $R$-admissible} for the pair $(T,{\bf k})$ as follows.
\begin{itemize}
\item For each vertex $r \in R$ of $T$, there exists a unique edge $i(r)'$ of $T$  incident with $r$ such that $k_r+r=l_{i(r)}+i(r)$.
\item For the edges $j'=(p,q)$ that do not appear in the image $i(R)'$ we have 
$\min(k_p+p,k_q+q) \le l_j+j < \max(k_p+p,k_q +q)$. (Note that 
for those edges we do not allow $l_j+j=k_p+p$ or $l_j+j=k_q+q$ if $p \in R$ or 
$q \in R$, respectively.)
\end{itemize}
The vector $\bf l$ is said to be {\bf $R$-admissible} if the function $i:R \to [n-1]$ is injective. If the function is not injective then we choose for each pair of distinct vertices $r, s \in R$ that share an edge $i(r)'=i(s)'$ one endpoint to be the dominating endpoint. 
\end{defi}

\medskip

\begin{figure}
\begin{center}
\mbox{\scalebox{0.40}{%
    \includegraphics{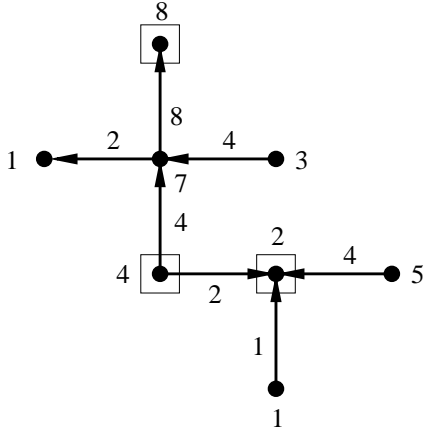}}}
\end{center}
\caption{\label{rad} An $R$-admissible labeling; the vertices of $R$ are enclosed by squares.}   
\end{figure}

An example is given in Figure~\ref{rad}. 
For the extreme cases concerning $R$, we have the following: the weak  $\emptyset$-admissibility coincides
with the ordinary admissibility and there exists no 
$[n]$-admissible vector as there is no injective function $i:[n] \to [n-1]$. If
$n=1$ then there exists an $R$-admissible vector if and only if $R=\emptyset$, namely the empty set.

\medskip

We introduce the {\bf sign} which we  associate with $(T,{\bf k})$, $i:R \to [n-1]$ and a choice of dominating vertices (if necessary). The following manner of speaking will turn out to be useful: if we refer to the minimum of an edge then we mean the minimum of the two labels of the endpoints of the edge or, by abuse of language, the respective vertex where this minimum is attained; similar for 
the maximum. If, for an edge $j'$, the labels on the two 
endpoints coincide then the edge must be in the image $i(R)'$. If $i^{-1}(j)$ contains a unique vertex then 
we define this to be the ``maximum'' of the edge and if $i^{-1}(j)$ contains both endpoints then the dominating vertex is defined as the ``maximum''; in both cases the other endpoint is defined as the minimum. As for the sign, we let each vertex that is an inversion contribute a $-1$ (which is the case when it is directed from its maximum to its minimum) as well as each $r \in R$ that is the minimum of the edge $i(r)'$.

\medskip

We define {\bf $(n,m,R)$-Gelfand-Tsetlin tree sequence} as follows.

\begin{defi}
\label{seq}
 Let $m \le n$ be positive integers, 
${\mathcal T}=(T_1,\ldots,T_n)$ be a tree sequence, $R \subseteq [m]$ be 
a set of vertices of $T_m$ and ${\bf k} \in \mathbb{Z}^n$ be a shifted labeling of the tree $T_n$.  An {\bf $(n,m,R)$-Gelfand-Tsetlin tree sequence} associated with ${\mathcal T}$ and 
${\bf k}$ is a sequence
${\bf L} = ({\bf l}_1,{\bf l}_2,\ldots,{\bf l}_n)$ with ${\bf l}_i \in \mathbb{Z}^i$ and 
${\bf l}_n={\bf k}$ which has  the following properties. 
\begin{itemize}
\item  The shifted labeling ${\bf l}_{i-1}$ is admissible for the pair $(T_i,{\bf l}_i)$ if $i \in \{2,3,\ldots,n\} \setminus \{m\}$.
\item The shifted labeling ${\bf l}_{m-1}$ is weakly $R$-admissible for the pair 
$(T_m,{\bf l}_m)$. 
\end{itemize}
If the function $i:R \to [m-1]$, which manifests the weak $R$-admissibility 
is not injective then the {\bf $(n,m,R)$-Gelfand-Tsetlin tree sequence} comes 
along with a set of dominating vertices as described in Definition~\ref{weakad}; all choices are possible.
We let ${\mathcal L}_{n,m,R}({\mathcal T},{\bf k})$ denote the set of these sequences. For an integer $\rho \le 
m$, we denote by  ${\mathcal L}_{n,m,\rho}({\mathcal T},{\bf k})$ the union over all
$\rho$-subsets $R$ of $[m]$.
Concerning 
the sign, we define
$$
\sgn {\bf L} = (-1)^\text{$\#$ of inversion of $\bf L$} \cdot 
 (-1)^\text{$\#$ of vertices $r \in R$ s.t. $r$ is the minimum of $i(r)'$} \cdot \sgn {\mathcal T}.
$$
We let $L_{n,m,R}({\mathcal T},{\bf k})$, respectively $L_{n,m,\rho}({\mathcal T},{\bf k})$,  denote the signed enumeration 
of these objects.
\end{defi}

The following is obvious but crucial: the quantity $L_{n,m,R}({\mathcal T},{\bf k})$ does not change if we pass in Definition~\ref{seq} from 
weak $R$-admissibility to $R$-admissibility as changing the dominating vertex from one endpoint of a shared edge to the other is a sign-reversing involution.

\medskip

We are in the position to give the combinatorial interpretation for the expression in \eqref{subset}.  In order to state the result, we introduce a convenient notation: if $R=\{i_1,\ldots,i_{\rho}\} \subseteq [n]$ then
$$
\Delta_{{\bf k}_R} f({\bf k}): = \Delta_{k_{i_1}} \cdots \Delta_{k_{i_{\rho}}} f(k_1,\ldots,k_n).$$
(The analog convention for $E_{{\bf k}_R} f({\bf k})$ will be used below.)

\begin{prop} 
\label{first}
Let  $R \subseteq [n]$. Then 
$\Delta_{{\bf k}_R} {L}_n({\mathcal T},{\bf k}) = L_{n,n,R}({\mathcal T},{\bf k})$.
\end{prop}

\medskip

This immediately implies the following combinatorial interpretation for the left-hand side of \eqref{toshow}.

\begin{cor} 
\label{firstcor}
Let $\rho \in \{0,1,\ldots,n\}$. Then
$e_{\rho}(\Delta_{k_1},\ldots,\Delta_{k_n}) {L}_n({\mathcal T},{\bf k}) = L_{n,n,\rho}({\mathcal T},{\bf k}).$
\end{cor}

\medskip

The following lemma is used in several places of our proofs in the remainder of this section. 

\begin{lem} 
\label{distinct}
For an integer $t < m$, we fix a set $P$ of pairs of edges of $T_{t+1}$ and
let ${\mathcal L}_{n,m,R,P}({\mathcal T}, {\bf k})$ denote the subset of labeled tree sequences in
${\mathcal L}_{n,m,R}({\mathcal T}, {\bf k})$ such that for each pair in 
$P$ the edge labels of the respective edges of $T_{t+1}$ are distinct. 
Then the signed 
enumeration of this subset is equal to the signed enumeration of the whole set.
\end{lem} 

\medskip

{\it Proof.} We consider the complement 
of ${\mathcal L}_{n,m,R,P}({\mathcal T}, {\bf k})$ and suppose that for $(i,j) \in P$ the edge labeling
${\bf l}_{t}+(1,2,\ldots,t)$ of $T_{t+1}$ is equal in the coordinates
$i$ and $j$. If there is more than one pair then we choose the pair which is minimial with respect to the fixed order on $P$. Then, we may replace the tree $T_t$ in ${\mathcal T}$ by a tree where vertex $i$ and $j$ are adjacent. The assertion follows as such as tree does not possess an admissible edge labeling. \qed

\medskip

{\it Proof of Proposition~\ref{first}.}  We consider subsets of ${\mathcal L}_{n}({\mathcal T},{\bf k})$ indexed by two disjoint subsets $P, Q  \subseteq [n]$ of vertices of $T_n$: let ${\mathcal L}_{n}({\mathcal T},{\bf k},P,Q)$ denote the set of Gelfand-Tsetlin tree sequences in  ${\mathcal L}_{n}({\mathcal T},{\bf k})$ such that
for the edge labeling ${\bf l} \in \mathbb{Z}^{n-1}$ of the largest tree $T_n$ in the tree sequence $\mathcal T$ the following is fulfilled:
\begin{itemize}
\item For each $p \in P$, there exists an edge $i(p)'$ of $T_n$ incident with $p$ such 
that $k_p+p$ is the minimum of $i(p)'$ and
$l_{i(p)}+i(p)=k_p+p$. 
\item For each $q \in Q$, there exists an edge $i(q)'$ in $T_n$ incident with $q$ such that $k_q+q$ is the maximum of $i(q)'$ and $l_{i(q)}+i(q)=k_q+q-1$.
\end{itemize}
We denote the respective signed enumeration by ${L}_{n}({\mathcal T},{\bf k},P,Q)$. Suppose $r \notin P, Q$. Then 
$$
\Delta_{k_r} {L}_{n}({\mathcal T},{\bf k},P,Q) = 
E_{k_r} {L}_{n}({\mathcal T},{\bf k},P,Q \cup \{r\})
-  {L}_{n}({\mathcal T},{\bf k},P \cup \{r\},Q).
$$
In order to see this, consider an element of $E_{k_r} {\mathcal L}_{n}({\mathcal T},{\bf k},P,Q)$ with the 
following property: for each edge $i'$ of $T_n$ that is incident with vertex $r$ of $T_n$ and such that the 
vertex label of the other endpoint of $i'$ is smaller than $k_r+r+1$ we have that the respective
edge label $l_i+i$ is smaller than $k_r+r$. In this case, we may change the vertex label of $r$ to $k_r+r$ to obtain an element of 
${\mathcal L}_{n}({\mathcal T},{\bf k},P,Q) \setminus {\mathcal L}_{n}({\mathcal T},{\bf k},P \cup \{r\},Q)$. Thus, these elements cancel in the difference on the left-hand side and we are left with the elements on the right-hand side.

\medskip

This implies by induction with respect to the size of $R \subseteq [n]$ that
\begin{equation}
\label{alternatingsum}
\Delta_{{\bf k}_R} 
{L}_{n}({\mathcal T},{\bf k}) = \sum_{Q \subseteq R}
(-1)^{|R|+|Q|} E_{{\bf k}_{Q}} {L}_{n}({\mathcal T},{\bf k},R \setminus Q,Q).
\end{equation}
The right-hand side is in fact equal to the signed enumeration of 
${\mathcal L}_{n,n,R}({\mathcal T},{\bf k})$: in order to see this, we may assume by 
Lemma~\ref{distinct} that the edge labels of $T_{n}$ are distinct, both 
in ${\mathcal L}_{n,n,R}({\mathcal T},{\bf k})$ and in $E_{{\bf k}_{Q}} {\mathcal L}_{n}({\mathcal T},{\bf k},R \setminus Q,Q)$. This implies that for each tree sequences in $E_{{\bf k}_{Q}} {\mathcal L}_{n}({\mathcal T},{\bf k},R \setminus Q,Q)$ and each $r \in R$, there is a unique edge $i(r)'$ of $T_n$ with $l_{i(r)}+i(r) = k_r + r$. Now, we may convert elements of $E_{{\bf k}_{Q}} {\mathcal L}_{n}({\mathcal T},{\bf k},R \setminus Q,Q)$ into elements of 
${\mathcal L}_{n,n,R}({\mathcal T},{\bf k})$ by decreasing the 
labels of the vertices in $Q$ by $1$. 
We obtain elements, where for $r \in Q$, 
the vertex label $k_r+r$ is the maximum of $i(r)'$ and, for $r \in R \setminus Q$, the vertex label 
$k_r+r$ is the minimum of $i(r)'$ -- attached 
with a sign according to the number cases where $k_r+r$ is the minimum of the edge $i(r)'$. The fact that the edge labels are distinct and since there always exists an edge label that is equal to $k_r+r$ implies that it is irrelevant that the intervals for the possible labels of the edges incident with $r$ were slightly changed when passing from $E_{{\bf k}_{Q}} {\mathcal L}_{n}({\mathcal T},{\bf k},R \setminus Q,Q)$ to ${\mathcal L}_{n,n,R}({\mathcal T},{\bf k})$.

\medskip

However, by decreasing the vertex label of a vertex $q \in Q$ of 
an element in $E_{{\bf k}_{Q}} {\mathcal L}_{n}({\mathcal T},{\bf k},R \setminus Q,Q)$ by $1$ to $k_q+q$, this value may reach the vertex label $k_p+p$ of a vertex $p$ that is adjacent to $q$; in this case we have to guarantee that $k_q+q$ can still be  identified as the maximum of the edge $j'$ connecting $p$ and $q$. The  assumption implies $i(q)=j$. 
If  $p \notin R$ then, when considering the labeled tree sequence as an element of ${\mathcal L}_{n,n,R}({\mathcal T},{\bf k})$,
 the vertex $q$ is the maximum of $j'$ by definition. If, on the other hand, $p \in R$, then we also have $i(p)=j$
 and we let $q$ be the dominating vertex of the edge to remember that it used to be the maximum of the edge $j'$. Thus it is clear how to reverse the procedure.
\qed

\medskip

In the definition of the $R$-admissibility, we have fixed a set $R$ of vertices of $T$. However, we may as well fix  the image $i(R)=:R'$ of the injective function $i:R \to [n-1]$, which corresponds to a set of edges of $T$. 

\begin{defi}
Let $T$ be an $n$-tree, ${\bf k} \in \mathbb{Z}^{n}$ and $R' \subseteq [n-1]$.  A vector ${\bf l} \in 
\mathbb{Z}^{n-1}$ together with a function $t:R' \to [n]$ is said to be
{\bf $R'$-edge-admissible} for the pair $(T,{\bf k})$ if 
${\bf l}$ is $t(R')$-admissible for the pair $(T,{\bf k})$, where $t^{-1}:t(R') \to [n-1]$ is the function that proves the $t(R')$-admissibility.
\end{defi}

In analogy to Definition~\ref{seq}, it is also clear how to define Gelfand-Tsetlin tree sequences associated with a triple $(n,m,R')$, where $m \le n$ are positive integers and $R' \subseteq [m-1]$ corresponds to a subset of {\it edges} of $T_m$. We denote this set by ${\mathcal L}_{n,m}^{R'}({\mathcal T},{\bf k})$ and by 
${L}_{n,m}^{R'}({\mathcal T},{\bf k})$ its signed enumeration. Note that ${\mathcal L}_{n,m,\rho}({\mathcal T},{\bf k})$ is also the union of  ${\mathcal L}_{n,m}^{R'}({\mathcal T},{\bf k})$, where $R'$ is a 
$\rho$-subset of $[m-1]$.  

\medskip

In the proof of the next proposition, it will be helpful to replace the $R'$-edge-admissibility in the definition of ${L}_{n,m}^{R'}({\mathcal T},{\bf k})$
by a more general notion, which we call 
{\it weak $R'$-edge-admissibility} and define as follows.

\medskip

\begin{defi}
Let $T$ be an $n$-tree, ${\bf k} \in \mathbb{Z}^{n}$ and $R' \subseteq [n-1]$.  A vector ${\bf l} \in 
\mathbb{Z}^{n-1}$  is said to be
{\bf weakly $R'$-edge-admissible} for the pair $(T,{\bf k})$ if there exists a function $t:R' \to [n]$ such that the following conditions are fulfilled.
\begin{itemize}
\item For all $r \in R$, the edge $r'$ of $T$ is incident with the vertex $t(r)$ of $T$ and 
$l_r+r=k_{t(r)}+t(r)$.
\item For all $r \in [n-1] \setminus R'$, we have $\min(k_p+p,k_q+q) \le l_r +r < \max(k_p+p,k_q+q)$, where $r'=(p,q)$ in T.
\end{itemize} 
The sign we associate is defined as follows: each inversion contributes 
a $-1$ as well as each edge $r'$ of $R'$ such that $t(r)$ is the minimum of the edge. (If the two vertex labels of an edge $r'$ coincide then  it must be an 
element of $R'$ and we define $t(r)$ as the ``maximum'' of the edge.)
\end{defi}

\medskip

To obtain the ordinary edge-admissibility we have to require in addition that for all $r \in R'$ the following is fulfilled: suppose $s'$ is an edge of $T$ incident with vertex $t(r)$ such that $l_s+s=k_{t(r)}+t(r)$ then we have $r=s$. However, 
the violation of this condition would require two edges of $T$ to have the same label, which can be avoided for an element of ${L}_{n,m}^{R'}({\mathcal T},{\bf k})$ by the argument given in Lemma~\ref{distinct}. 
\medskip

The following proposition will finally imply \eqref{toshow}.

\begin{prop} 
\label{second}
Let $R \subseteq [m-1]$. Then 
$L_{n,m}^{R}({\mathcal T},{\bf k}) = L_{n,m-1,R}({\mathcal T},{\bf k})$.
\end{prop}

An immediate consequence is the following.

\begin{cor} Let $\rho$ be a non-negative integer. Then 
$L_{n,m,\rho}({\mathcal T},{\bf k}) = L_{n,m-1,\rho}({\mathcal T},{\bf k})$.
\end{cor}

The corollary implies 
$
L_{n,m,\rho}({\mathcal T},{\bf k}) = 0
$
if $\rho$ is non-zero as  ${\mathcal L}_{n,m,\rho}({\mathcal T},{\bf k})=\emptyset$ if $\rho \ge m$, since there is no injective function from $[\rho]$ to $[m-1]$. 
By Corollary~\ref{firstcor}, \eqref{toshow} finally follows.

\medskip

{\it Proof of Proposition~\ref{second}.} We restrict our considerations to the case that $m=n$ as the general case is analog.  By 
 Lemma~\ref{distinct}, we assume that the edge labels of $T_{n-1}$ are distinct, both in 
 ${\mathcal L}_{n,n}^{R}({\mathcal T},{\bf k})$ and in 
 ${\mathcal L}_{n,n-1,R}({\mathcal T},{\bf k})$.

\medskip

\begin{figure}
\begin{center}
\mbox{\scalebox{0.30}{%
    \includegraphics{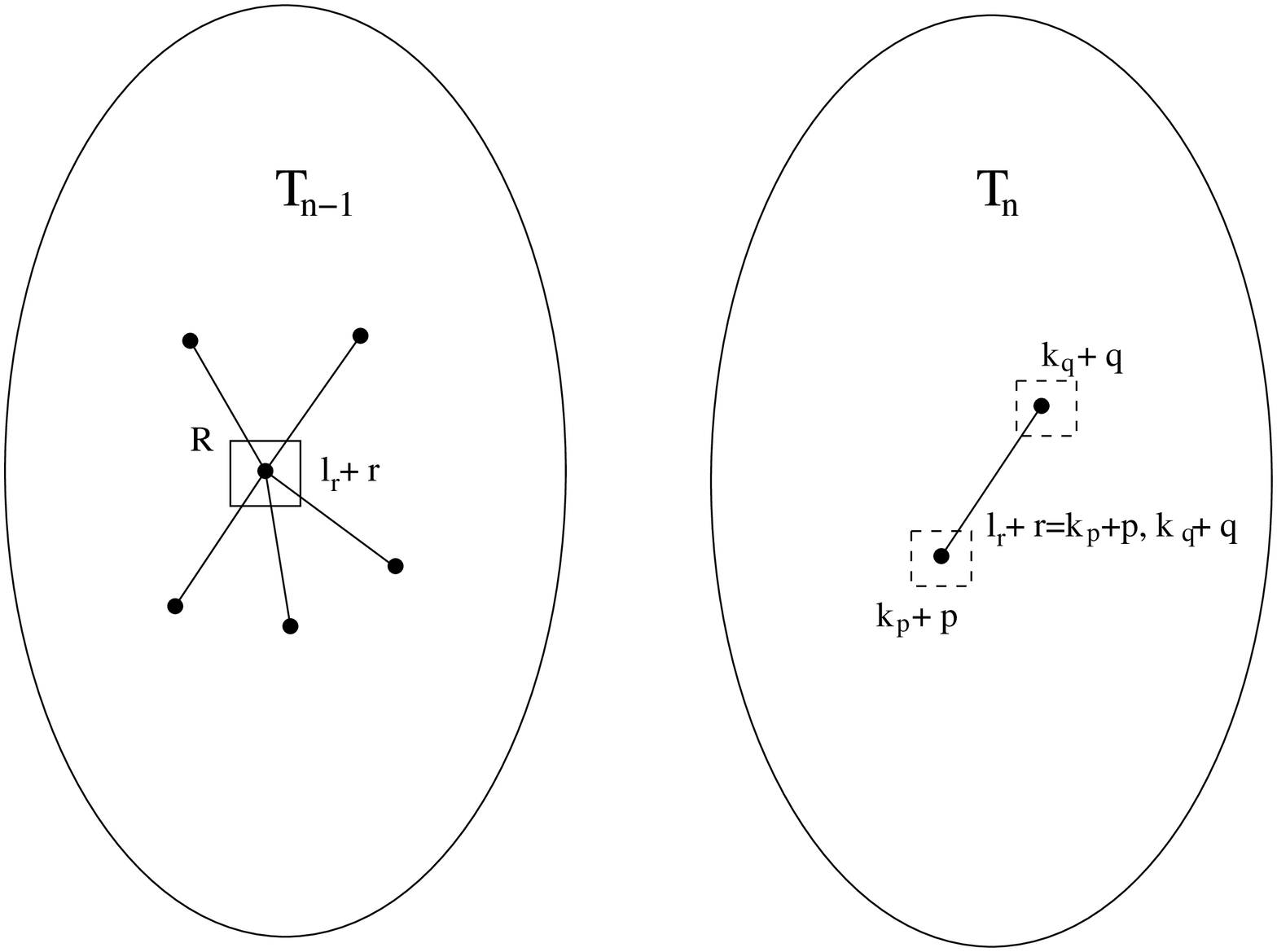}}}
\end{center}
\caption{\label{skizze} Situation in the proof of Proposition~\ref{second}.}   
\end{figure}

We consider an element 
of ${\mathcal L}_{n,n}^{R}({\mathcal T},{\bf k})$,  
denote by ${\bf l} \in \mathbb{Z}^{n-1}$ the respective shifted edge labeling of $T_{n}$ and by
$t:R \to [n]$ the function that proves the weak $R$-edge-admissibility of the vector $\bf l$ for the pair 
$(T_n,{\bf k})$. Suppose $r \in R$ and that $p,q$ are the vertices of the edge $r'$ in $T_n$  then  
we have either $t(r)=p$ or $t(r)=q$. We denote the first subset of 
${\mathcal L}_{n,n}^{R}({\mathcal T},{\bf k})$ by $M_{r,p}$ and the second subset 
by $M_{r,q}$. The situation is sketched in Figure~\ref{skizze}.

\medskip

Assuming w.l.o.g. that $k_p+p \le k_q +q$, we first observe that we can restrict our attention to the case that there is at least on edge incident with vertex $r$ in $T_{n-1}$, the label of which lies in the interval $[k_p+p,k_{q}+q)$. This is because for the other elements, $l_r+r \to k_q+q$ and $t(r) \to q$ induces a sign reversing bijection from $M_{r,p}$ to $M_{r,q}$. In the following, we address these edges as the relevant edges of $r$.

\medskip

In order to construct an element of  ${\mathcal L}_{n,n-1,R}({\mathcal T},{\bf k})$ we perform the following shifts to the  labels $l_r+r$ for all $r \in R$: if the fixed element of ${\mathcal L}_{n,n}^{R}({\mathcal T},{\bf k})$ is an element of $M_{r,p}$, we shift $l_{r}+r$ to the minimum of incident edge labels in $T_{n-1}$ no smaller than $k_p+p$
and let $j(r)'$
be the respective edge, while for elements of $M_{r,q}$ we shift $l_{r}+r$ to the maximum of incident edge labels in $T_{n-1}$ smaller than $k_q+q$ and let $j(r)'$ be the respective edge. These edges $j(r)$ are unique as the edge labels are assumed to be distinct. The contribution of $-1$ to the sign of the elements in $M_{r,p}$ that comes from
the fact that the edge label of $r'$ in $T_n$ is equal to the minimum of the edge translates in the new element into the contribution of $-1$ of the edge $j(r)'$ in $T_{n-1}$ as 
its edge label is  also equal to the minimum of the edge. If this procedure causes two distinct vertices $r, s \in R$ to share an edge $j(r)'=j(s)'$ then we let the dominating vertex  be the maximum of the respective edge in the 
original element.
\medskip

The precise description of the elements in ${\mathcal L}_{n,n-1,R}({\mathcal T},{\bf k})$ that appear as a result of 
this procedure is the following. For each $r \in R$, one of the following two possibilities applies: suppose $p,q$ are the endpoints of $r'$ in $T_n$ and w.l.o.g. $k_p+p \le k_q + q$ then either
\begin{itemize}
\item the vertex $r$ is the minimum of the edge $j(r)'$ and the edge label of $j(r)'$ 
is the minimum under all relevant edges of $r$, or 
\item the vertex $r$ is the maximum of the edge $j(r)'$ and 
the edge label of $j(r)'$ 
is the 
maximum under all relevant edges of $r$.
\end{itemize}
For such an element it is also clear how to invert the procedure to reobtain 
an element of ${\mathcal L}_{n,n}^{R}({\mathcal T},{\bf k})$.

\medskip

Finally, we define a sign-reversing involution on the set of elements of ${\mathcal L}_{n,n-1,R}({\mathcal T},{\bf k})$
that do not fulfill this requirement: suppose that $r \in R$ is minimal such that the requirement is not met and that $r$ is the minimum of the edge $j(r)'$. 
Let $i'$ be the relevant edge of $r$, the edge label of which is maximal with the property that it is smaller than $l_r+r$. We shift $l_r+r$ to this edge 
label and set $j(r)=i$. If necessary we choose the dominating vertices such that the set of inversions remains unaffected. Then, $r$ is the maximum of the edge $j(r)'$. Likewise when $r$ is the maximum of the edge. The fact that we only work with relevant edges guarantees that we are able to perform the shift accordingly for the edge label of $r'$ in $T_n$.
\qed

\medskip

To conclude this section, we demonstrate that also \eqref{toshow} implies that ${L}_n({\mathcal T},{\bf k})$ is a polynomial in $k_1,\ldots, k_n$ of degree no greater than $n-1$ in every $k_i$.

\begin{lem}
\label{degree}
Suppose that $A(k_1,\ldots,k_n)$ is a function with 
$$
e_{\rho}(\Delta_{k_1},\ldots,\Delta_{k_n}) A(k_1,\ldots,k_n) = 0
$$
for all $\rho >0$. Then $\Delta^n_{k_i} A(k_1,\ldots,k_n) = 0$ for all $i \in \{1,2,\ldots,n\}$.
\end{lem}

{\it Proof.} We define 
$$A_{\rho,i}(k_1,\ldots,k_n) = e_{\rho}(\Delta_{k_1},\ldots,\widehat{\Delta_{k_i}},\ldots,\Delta_{k_n}) A(k_1,\ldots,k_n),$$ 
where $\widehat{\Delta_{k_i}}$ indicates that $\Delta_{k_i}$ does not appear in the argument.
We use the identity 
$$
e_{\rho}(X_1,\ldots,X_n) = e_{\rho}(X_1,\ldots,\widehat{X_i},\ldots,X_n) + X_i e_{\rho-1}(X_1,\ldots,\widehat{X_i},\ldots,X_n)
$$
and the assumption to see that 
$$
A_{\rho,i}(k_1,\ldots,k_n) = - \Delta_{k_i} A_{\rho-1,i}(k_1,\ldots,k_n).
$$
This implies 
$$
A_{\rho,i}(k_1,\ldots,k_n) = (-1)^{\rho} \Delta^{\rho}_{k_i} A(k_1,\ldots,k_n)
$$
by induction with respect to $\rho$. As $A_{n,i}(k_1,\ldots,k_n)=0$, the assertion follows. \qed

\section{Monotone triangles}
\label{concluding}

I would like to see an analog ``theory'' for monotone triangles (Gelfand-Tsetlin patterns with strictly increasing rows), which seems conceivable as there are several properties of the unrestricted patterns for which we have a corresponding (though in some cases more complicated) property of monotone triangles. For instance, it is known \cite{monotonetriangles} that the number $\alpha(n;k_1,\ldots,k_n)$ of monotone 
triangles with bottom row $k_1,k_2,\ldots,k_n$ is given by
\begin{equation}
\label{operator}
\prod_{1 \le p < q \le n} (E_{k_p} + E_{k_q}^{-1} - E_{k_p} E^{-1}_{k_q}) 
\prod_{1 \le i < j \le n} \frac{k_j - k_i + j -i}{j-i} = 
\prod_{1 \le p < q \le n} (\id + \Delta_{k_p} \delta_{k_q})
\prod_{1 \le i < j \le n} \frac{k_j - k_i + j -i}{j-i},
\end{equation}
where $\delta_x := \id - E^{-1}_{x}$. To start with, we give four different  (but related)
combinatorial extensions of  $\alpha(n;k_1,\ldots,k_n)$ to all $(k_1,\ldots,k_n) \in \mathbb{Z}^n$ in this section, and 
then present certain other properties of  $\alpha(n;k_1,\ldots,k_n)$, for which it would be nice to have 
combinatorial proofs of the type as we have presented them in this article for Gelfand-Tsetlin patterns. This is because these properties imply, on the one hand, \eqref{operator} and, on the other hand, the refined alternating sign matrix theorem. The latter will be explained at the end of this section.

\medskip

\subsection{Four combinatorial extensions of $\alpha(n;k_1,\ldots,k_n)$ to all $(k_1,\ldots,k_n) \in \mathbb{Z}^n$}
The quantity $\alpha(n;k_1,\ldots,k_n)$ obviously satisfies the following recursion for any sequence $(k_1,k_2,\ldots,k_n)$  of strictly increasing integers.
\begin{equation}
\label{recursion}
\alpha(n;k_1,\ldots,k_n) = 
\sum_{(l_1,\ldots,l_{n-1}) \in \mathbb{Z}^{n-1} \atop 
k_1 \le l_1 \le k_2 \le l_2 \le k_3 \le \ldots \le k_{n-1} \le l_{n-1} \le k_n, l_i \not= l_{i+1}} 
\alpha(n-1;l_1,\ldots,l_{n-1})
\end{equation}
To obtain an extension of the combinatorial interpretation of $\alpha(n;k_1,\ldots,k_n)$ to all 
$(k_1,\ldots,k_n) \in \mathbb{Z}^n$, it is convenient to write this summation in terms of ``simple'' summations 
$\sum\limits_{i=a}^b f(i)$, i.e. summations over intervals. This is because we can then use the extended definition of the summation, i.e.
$\sum\limits_{i=a}^{a-1} f(i)=0$ and $\sum\limits_{i=a}^b f(i)=-\sum\limits_{i=b+1}^{a-1} f(i)$ if $b+1 \le a-1$. Note that if $p(i)$ is a polynomial in $i$ then there exists a polynomial $q(i)$ with $\Delta_i q(i) = p(i)$, which implies $\sum\limits_{i=a}^{b} p(i) = q(b+1) - q(a)$ if $a \le b$ and, consequently, that this sum is a polynomial in $a$ and $b$. The extension of the simple summation we have just introduced was chosen such that the latter identity is true for all $a,b \in \mathbb{Z}$.  After we have given at least one representation of the summation in \eqref{recursion} in terms of simple summations, this shows  that $\alpha(n;k_1,\ldots,k_n)$ can be represented by a polynomial in $k_1, k_2, \ldots, k_n$ if $k_1 < k_2 < \ldots < k_n$. (This polynomial is in fact unique as a polynomial in $k_1,k_2,\ldots,k_n$ is uniquely determined by its values on the set of $n$-tuples $(k_1,k_2,\ldots,k_n) \in \mathbb{Z}^n$ with $k_1 < k_2 < \ldots < k_n$.) The extended monotone triangles with prescribed bottom row $k_1,k_2,\ldots,k_n$ will be chosen such that these objects are enumerated by this polynomial for all $(k_1,\ldots,k_n) \in \mathbb{Z}^n$. In particular, it will certainly not be the naive extension, which sets $\alpha(n;k_1,\ldots,k_n)=0$ if $k_1,k_2,\ldots,k_n$ is not strictly increasing. 

\medskip

\subsubsection{First extension} If we assume that $k_1 < k_2 < \ldots < k_n$, then one possibility to write the summation in \eqref{recursion} in terms of simple summations is the following: we choose a 
subset $\{l_{i_1},l_{i_2},\ldots,l_{i_p} \} \subseteq \{l_1,\ldots,l_{n-1}\}$ for which we have $l_{i_j}=k_{i_j}$. For all other $l_q$ we have $k_{q} < l_q \le k_{q+1}$, except for the case that $q+1=i_j$ for a $j$, where we have $k_{q} < l_q < k_{q+1}$.  More formally, 
$$
\sum_{p \ge 0} \sum_{1 \le i_1 < i_2 < \ldots < i_p \le n-1} 
\sum_{l_1=k_1+1}^{k_2} \sum_{l_2=k_2+1}^{k_3} \ldots 
\sum_{l_{i_1-1}=k_{i_1-1}+1}^{k_{i_1}-1} \sum_{l_{i_1}=k_{i_1}}^{k_{i_1}} \ldots
\sum_{l_{i_p-1}=k_{i_p-1}+1}^{k_{i_p}-1} \sum_{l_{i_p}=k_{i_p}}^{k_{i_p}} \ldots
\sum_{l_{n-1}=k_{n-1}+1}^{k_n}
$$
where in the exceptional case that $i_j=i_{j-1}+1$ the expression 
$\sum\limits_{l_{i_j-1}=k_{i_j-1}+1}^{k_{i_j}-1} \sum\limits_{l_{i_j}=k_{i_j}}^{k_{i_j}}$ is replaced by
$\sum\limits_{l_{i_{j-1}}=k_{i_{j-1}}}^{k_{i_{j-1}}} \sum\limits_{l_{i_j}=k_{i_j}}^{k_{i_j}}$.
This leads to the following extension: a monotone triangle of order $n$ is a triangular array $(a_{i,j})_{1 \le j \le i \le n}$ of integers such that the following conditions are fulfilled.
\begin{itemize}
\item There is a subset of special entries $a_{i,j}$ with $i < n$ for which we require $a_{i,j}=a_{i+1,j}$. We mark these entries with a star on the left.
\item If $a_{i,j}$ is not a special entry then we have to distinguish between the case that $a_{i,j}$ is the left neighbour of a special entry or not.
\begin{itemize}
\item If $a_{i,j+1}$ is not special (which includes also the case that $a_{i,j+1}$ does not exist) then 
$a_{i+1,j} < a_{i,j} \le a_{i+1,j+1}$ in case that $a_{i+1,j} < a_{i+1,j+1}$ and $a_{i+1,j+1} < a_{i,j} \le a_{i+1,j}$ otherwise. (There exists no pattern with $a_{i+1,j}=a_{i+1,j+1}$.) In the latter case we have an inversion.
\item If $a_{i,j+1}$ is special then $a_{i+1,j} < a_{i,j} <  a_{i+1,j+1}$ or 
$a_{i+1,j+1} \le a_{i,j} \le a_{i+1,j}$. (There exists no pattern with $a_{i+1,j+1}=a_{i+1,j}+1$.) In the latter case we have an inversion.
\end{itemize}
\end{itemize}
The sign of a monotone triangle is $-1$ to the number of inversions. Then $\alpha(n;k_1,\ldots,k_n)$ is the signed enumeration of monotone triangles with $a_{n,i}=k_i$.  Here is an example of such an array.
\begin{center}
\begin{tabular}{ccccccccccccccccc}
  &   &   &   &   &   &   &   & $3$ &   &   &   &   &   &   &   & \\
  &   &   &   &   &   &   & $^*2$ &   & $6$ &   &   &   &   &   &   & \\
  &   &   &   &   &   & $2$ &   & $4$ &   & $^*6$ &   &   &   &   &   & \\
  &   &   &   &   & $3$ &   & $^*1$ &   & $6$ &   & $7$ &   &   &   &   & \\
  &   &   &   & $3$ &   & $1$ &  &   $7$ &   & $5$   &  & $8$  &   &   &   & 
\end{tabular}
\end{center}

\medskip

\subsubsection{Second extension}
The summation can also be written in the following more symmetric manner: we choose a subset $I \subseteq [n-1]$ such that $l_i=k_i$ if $i \in I$ and a subset $J \subseteq [n-1]$ such that $l_j=k_{j+1}$ if $j \in J$. The sets $I, J$ have to be disjoint and, moreover, $i \in I$ implies $i-1 \notin J$ (which is equivalent to $(I-1) \cap J = \emptyset)$. On the other hand, if $h \in [n-1] \setminus (I \cup J)$ then $k_h < l_h < k_{h+1}$. Equivalently, 
$$
\sum_{p, q \ge 0} \sum_{I=\{i_1,\ldots,i_p\}, J=\{j_1,\ldots,j_q\} \subseteq [n-1] \atop I \cap J = \emptyset, (I-1) \cap J = \emptyset}
\sum_{l_{i_1}=k_{i_1}}^{k_{i_1}} \ldots \sum_{l_{i_p}=k_{i_p}}^{k_{i_p}} \sum_{l_{j_1}=k_{j_1+1}}^{k_{j_1+1}} \ldots \sum_{l_{j_q}=k_{j_q+1}}^{k_{j_q+1}}
\sum_{l_{h_1}=k_{h_1}+1}^{k_{h_1+1}-1} \ldots \sum_{l_{h_r}=k_{h_r}+1}^{k_{h_r+1}-1},
$$
where $[n-1] \setminus (I \cup J) = \{h_1,\ldots,h_r\}$. Using this representation, we can deduce the following extension: a monotone triangle of order $n$ is a triangular array $(a_{i,j})_{1 \le j \le i \le n}$ of integers such that the following conditions are fulfilled.
\begin{itemize}
\item There is a subset of ``left-special'' entries $a_{i,j}$ with $i < n$ for which we require $a_{i,j}=a_{i+1,j}$ and we mark them 
with a star on the left as well as a subset of ``right-special'' entries $a_{i,j}$ with $i <n$ for which we require $a_{i,j}=a_{i+1,j+1}$ and mark them with a star on the right. 
\item An entry can not be a left-special entry and a right-special entry. If a right-special entry and a left-special entry happen to be in the same row then the right-special entry may not be situated immediately to the left of the left-special entry.
\item If $a_{i,j}$ is not a special entry then we have $a_{i+1,j} < a_{i,j} <  a_{i+1,j+1}$ or 
$a_{i+1,j+1} \le a_{i,j} \le a_{i+1,j}$, respectively. In the latter case we have an inversion.
\end{itemize} 
Next we give an example of such an array.
\begin{center}
\begin{tabular}{ccccccccccccccccc}
  &   &   &   &   &   &   &   & $3$ &   &   &   &   &   &   &   & \\
  &   &   &   &   &   &   & $2$ &   & $4$ &   &   &   &   &   &   & \\
  &   &   &   &   &   & $3$ &   & $^*2$ &   & $6$ &   &   &   &   &   & \\
  &   &   &   &   & $^*3$ &   & $2$ &   & $5^*$ &   & $7$ &   &   &   &   & \\
  &   &   &   & $3$ &   & $1$ &  &   $7$ &   & $5$   &  & $8$  &   &   &   & 
\end{tabular}
\end{center}
The sign of a monotone triangle is again $-1$ to the number of inversions and $\alpha(n;k_1,\ldots,k_n)$ is the signed enumeration of these extended monotone triangles with prescribed $a_{n,i}=k_i$.
Although we think that the fourth extension is probably the nicest, the first two extensions are the only ones where in case that $k_1 < k_2 < \ldots < k_n$ the removal of all stars leads to a monotone triangle in the original sense and no array is assigned a minus sign, i.e. we have a plain enumeration.

\medskip

\subsubsection{Third extension} Another possibility to write the summation in \eqref{recursion} in terms of simple summations  is the following. 
$$
\sum_{p \ge 0} (-1)^p \sum_{2 \le i_1 < i_2 < \ldots < i_p \le n-1 \atop i_{j+1} \not= i_j+1} 
\sum_{l_1=k_1}^{k_2} \sum_{l_2=k_2}^{k_3} \ldots 
\sum_{l_{i_1-1}=k_{i_1}}^{k_{i_1}} \sum_{l_{i_1}=k_{i_1}}^{k_{i_1}} \ldots
\sum_{l_{i_p-1}=k_{i_p}}^{k_{i_p}} \sum_{l_{i_p}=k_{i_p}}^{k_{i_p}} \ldots
\sum_{l_{n-1}=k_{n-1}}^{k_n}
$$
This leads to the following extension: a monotone triangle of order $n$ is a triangular array $(a_{i,j})_{1 \le j \le i \le n}$ of integers such that the following conditions are fulfilled. The entries 
$a_{i-1,j-1}$ and $a_{i-1,j}$ are said to be the parents of $a_{i,j}$. 
\begin{itemize}
\item Among the entries $(a_{i,j})_{1 < j < i \le n}$ we may have special entries such that if two of them happen to be in the same row they must not be adjacent. We mark these entries with a star. For the parents of a special entry $a_{i,j}$ we have require $a_{i-1,j-1}=a_{i,j}=a_{i-1,j}$.
\item If $a_{i,j}$ is not the parent of a special entry then $a_{i+1,j} \le a_{i,j} \le a_{i+1,j+1}$ and $a_{i+1,j+1} < a_{i,j} < a_{i+1,j}$, respectively. In the latter case we have an inversion.
\end{itemize}
In this case, the sign of a monotone triangle is $-1$ to the number of inversions plus the number of special entries. Then $\alpha(n;k_1,\ldots,k_n)$ is the signed enumeration of monotone triangles with $a_{n,i}=k_i$. Next we give an example of such an array.
\begin{center}
\Large 
\begin{tabular}{ccccccccccccccccc}
  &   &   &   &   &   &   &   & $\phantom{*} \atop 4$ &   &   &   &   &   &   &   & \\
  &   &   &   &   &   &   & $\phantom{*} \atop 5$ &   & $\phantom{*} \atop 3$ &   &   &   &   &   &   & \\
  &   &   &   &   &   & $\phantom{*} \atop 5$ &   & $\phantom{*} \atop 5$ &   & $\phantom{*} \atop 2$ &   &   &   &   &   & \\
  &   &   &   &   & $\phantom{*} \atop 3$ &   & $* \atop 5$ &   & $\phantom{*} \atop 2$ &   & $\phantom{*} \atop 2$ &   &   &   &   & \\
  &   &   &   & $\phantom{*} \atop 4$ &   & $\phantom{*} \atop 1$ &  &   $\phantom{*} \atop 7$ &   & $* \atop 2$   &  & $\phantom{*} \atop 5$  &   &   &   & 
\end{tabular}
\end{center}

\medskip

This is the extension that has already appeared in \cite{simplified}. There we have indicated that the non-adjacency requirement for special entries can also be ignored: suppose that $(a_{i,j})_{1 \le j \le i \le n}$ is an array with the properties given above accept that we allow special entries to be adjacent: suppose $a_{i,j}$ and $a_{i,j+1}$ are two adjacent special entries such that $i+j$ is maximal with this property. 
Then we have 
$a_{i-1,j-1} = a_{i,j} = a_{i-1,j} = a_{i,j+1} = a_{i-1,j+1}$. This 
implies that $a_{i-2,j-1} = a_{i-1,j} = a_{i-2,j}$ whether or not $a_{i-1,j}$ is a special entry, which implies that changing the status of the entry $a_{i-1,j}$ is a sign-reversing involution.

\medskip 

\subsubsection{Fourth extension} In order to explain the representation of \eqref{recursion} in terms of simple summations which is used for the third extension, it is convenient to 
use the operator
$V_{x,y}: = E_x^{-1} + E_y - E^{-1}_x E_y$. Then
$$
\sum_{ 
k_1 \le l_1 \le k_2 \le l_2  \le \ldots \le k_{n-1} \le l_{n-1} \le k_n, \atop l_i \not= l_{i+1}}  
a(l_1,\ldots,l_{n-1}) = \left. V_{k_1,k'_1} V_{k_2,k'_2} \cdots 
V_{k_n,k'_n} \sum_{l_1=k'_1}^{k_2} \sum_{l_2=k'_2}^{k_3} \ldots \sum_{l_{n-1}=k'_{n-1}}^{k_n} 
a(l_1,\ldots,l_{n-1}) \right|_{k'_i=k_i},
$$
if $k_1 < k_2 < \ldots < k_n$ is strictly increasing. (Note that $V_{k_1,k'_1}$ as well as $V_{k_n,k'_n}$ can also be removed as the application of $V_{x,y}$ 
to a function which does not depend on $x$ {\it and} $y$ acts as the identity. In order to convince oneself that this is indeed a valid representation of the summation in \eqref{recursion}, one can use induction with respect to $n$ to transform it into the representation of 
the first extension.)
This leads to the following extension, which we think is the nicest: a monotone triangle of order $n$ is an integer array $(a_{i,j})_{1 \le j \le i \le n}$ 
together with 
a function $f$ which assigns to each $a_{i,j}$ an element of 
$\{\leftarrow, \rightarrow, \leftrightarrow\}$ such that the  following conditions  are fulfilled for any element $a_{i,j}$ with $i < n$: we have to 
distinguish cases depending on the assignment of the arrows to the elements 
$a_{i+1,j}$ and $a_{i+1,j+1}$.
\begin{enumerate}
\item $f(a_{i+1,j})=\leftarrow$, $f(a_{i+1,j+1})=\leftarrow, \leftrightarrow$ : $a_{i+1,j} \le a_{i,j} < a_{i+1,j+1}$  or 
$a_{i+1,j+1} \le a_{i,j} < a_{i+1,j}$
\item $f(a_{i+1,j})=\leftarrow$, $f(a_{i+1,j+1})=\rightarrow$ : $a_{i+1,j} \le a_{i,j} \le a_{i+1,j+1}$  or 
$a_{i+1,j+1} < a_{i,j} < a_{i+1,j}$
\item $f(a_{i+1,j})=\leftrightarrow, \rightarrow$, $f(a_{i+1,j+1})=\leftarrow, \leftrightarrow$ : $a_{i+1,j} < a_{i,j} < a_{i+1,j+1}$  or 
$a_{i+1,j+1} \le a_{i,j} \le a_{i+1,j}$ 
\item $f(a_{i+1,j})=\leftrightarrow, \rightarrow$, $f(a_{i+1,j+1})=\rightarrow$ : $a_{i+1,j} < a_{i,j} \le a_{i+1,j+1}$  or 
$a_{i+1,j+1} < a_{i,j} \le a_{i+1,j}$.
\end{enumerate}
In Case~1 and Case~4, there exists no pattern if $a_{i+1,j}=a_{i+1,j+1}$, in Case~2, we have no pattern if  $a_{i+1,j}=a_{i+1,j+1}+1$ and, in Case~3, there is no pattern if $a_{i+1,j+1}=a_{i+1,j}+1$. In each case, we say that $a_{i,j}$ is an inversion if the second possibility applies. We define the sign of 
a monotone triangle to be $-1$ to the number of inversions plus the number of elements that are assigned the element ``$\leftrightarrow$''.  Then 
$\alpha(n;k_1,\ldots,k_n)$ is the signed enumeration of monotone triangles $(a_{i,j})_{1 \le j \le i \le n}$ of order $n$ with 
$a_{n,i}=k_i$. Here is an example.
\begin{center}
\Large
\begin{tabular}{ccccccccccccccccc}
  &   &   &   &   &   &   &   & $\leftrightarrow \atop 5$ &   &   &   &   &   &   &   & \\
  &   &   &   &   &   &   & $\leftarrow \atop 5$ &   & $\rightarrow \atop 6$ &   &   &   &   &   &   & \\
  &   &   &   &   &   & $\leftarrow \atop 4$ &   & $\leftarrow \atop 6$ &   & $\leftrightarrow \atop 7$ &   &   &   &   &   & \\
  &   &   &   &   & $\leftarrow \atop 2$ &   & $\leftarrow \atop 6$ &   & $\rightarrow \atop 7$ &   & $\rightarrow \atop 5$ &   &   &   &   & \\
  &   &   &   & $\rightarrow \atop 3$ &   & $\leftrightarrow \atop 1$ &  &   $\rightarrow \atop 8$ &   & $\rightarrow \atop 5$   &  & $\rightarrow \atop 4$  &   &   &   & 
\end{tabular}
\end{center}
In order to see that this extension comes from the presentation given above, note that, when expanding
$$
V_{k_1,k'_1} V_{k_2,k'_2} \cdots V_{k_n,k'_n} = 
(E^{-1}_{k_1} + E_{k'_1} - E^{-1}_{k_1} E_{k'_1})
(E^{-1}_{k_2} + E_{k'_2} - E^{-1}_{k_2} E_{k'_2}) \cdots
(E^{-1}_{k_n} + E_{k'_n} - E^{-1}_{k_n} E_{k'_n}),
$$
the assignment of ``$\leftarrow$'' to the entry $k_i$ in the bottom row corresponds to choosing $E^{-1}_{k_i}$ from the operator $V_{k_i,k'_i}$, while the assignment of ``$\rightarrow$'' to $k_i$ corresponds to choosing $E_{k'_i}$ and the assignment of ``$\leftrightarrow$'' corresponds to choosing $E^{-1}_{k_i} E_{k'_i}$. 

\medskip

In all cases, the combinatorial extension of  $\alpha(n;k_1,\ldots,k_n)$ is, generally speaking, a signed enumeration, which reduces to a  plain enumeration in the first and in the second case if $k_1, k_2, \ldots, k_n$ is strictly increasing. This can be generalized as 
follows.

\begin{prop}
\label{boundary}
Suppose $k_1,k_2,\ldots, k_n$ is a weakly increasing sequence of integers then $\alpha(n;k_1,\ldots,k_n)$ is the number of Gelfand-Tsetlin patterns with 
prescribed bottom row $k_1,\ldots,k_n$ and where all other rows are strictly increasing.
\end{prop} 

{\it Proof.} In order to see this, we use the first extension. Suppose $k_j=k_{j+1}$ and $(a_{i,j})_{1 \le j \le i \le n}$ is a respective pattern. As 
$a_{n,j}=a_{n,j+1}$ it follows that $a_{n-1,j}$ equal to this quantity as well and at least one of $a_{n-1,j}$ and $a_{n-1,j+1}$ must be special. We can exclude the latter possibility by the following sign reversing involution on the extended monotone triangles where $a_{n-1,j+1}$ is special in such a situation: let $j$ be maximal with this property. Then, changing the status of $a_{n-1,j}$ (from special to not special or vice versa) is a sign reversing involution. Thus we can assume that 
$a_{n-1,j+1}$ is not special (and, consequently, $a_{n-1,j}$ must be special) whenever we have $a_{n,j}=a_{n,j+1}$. 

\medskip

This can be used to show that $\alpha(n;k_1,\ldots,k_n)=0$ if there are $p,q$ with $1 \le p < q \le n-1$ such that $k_{p}=k_{p+1}$, 
$k_{q}=k_{q+1}$ and $k_{j}+1=k_{j+1}$ for $p < j < q$, which is one special case of the statement: as $a_{n-1,p+1}$ can be assume not to be special (which already settles the case $q=p+1$)  we can deduce that $a_{n-1,p+2}$ is not special (otherwise we would have no choice for $a_{n-1,p+1}$) and, by iterating this argument, we can see that $a_{n-1,j}$ is not special for $p+1 \le j \le q-1$. This implies that $a_{n-1,p+1}=a_{n,p+2}, a_{n-1,p+2}=a_{n,p+3}, \ldots, a_{n-1,q-1}=a_{n,q}$.
On the other hand, the fact that $a_{n-1,q}$ is special implies $a_{n-1,q-1}=a_{n,q-1}$, which is a contradiction.

\medskip

Thus we may assume that such $p,q$ do not exist for our sequence $k_1,k_2,\ldots,k_n$. Consequently, if $k_j=k_{j+1}$ then $k_{j-1} < k_j$ and 
$k_{j+1} < k_{j+2}$. As $a_{n-1,j}$ is  special and $a_{n-1,j+1}$ is not, we have $a_{n-1,j-1} < a_{n-1,j} < a_{n-1,j+1}$. \qed

\medskip

It should be remarked that the signed enumeration in the first and in the second extension is in general not a plain enumeration if $k_1, \ldots, k_n$ is 
weakly increasing but not strictly increasing. Also note that the proposition is equivalent to the fact that, for weakly increasing sequences $k_1,k_2,\ldots,k_n$, the application of the summation in \eqref{recursion} to $\alpha(n-1;l_1,\ldots,l_{n-1})$ is equivalent to the application of the representation of this summation in terms of simple summations  to $\alpha(n-1;l_1,\ldots,l_{n-1})$. (If the sequence is not increasing then the summation in \eqref{recursion} is over the emptyset and therefore zero.) As a next step, it would be interesting to figure out whether there is a notion analog to that of Gelfand-Tsetlin tree sequences for monotone triangles. This could be helpful in understanding the properties of $\alpha(n;k_1,\ldots,k_n)$, which we list next.

\subsubsection{Properties of $\alpha(n;k_1,\ldots,k_n)$}
In previous papers we have shown that $\alpha(n;k_1,\ldots,k_n)$ has the following properties.
\begin{enumerate}
\item For $n \ge 1$ and $i \in \{1,2,\ldots, n-1\}$, we have 
$$(\id + E_{k_{i+1}} E^{-1}_{k_i} S_{k_i,k_{i+1}}) V_{k_i,k_{i+1}} 
\alpha(n;k_1,\ldots,k_n) = 0.$$
(This is proved in \cite{monotonetriangles}.) 
\item For $n \ge 1$ and $i \in \{1,2,\ldots,n\}$, we have $\deg_{k_i} \alpha(n;k_1,\ldots,k_n) \le n-1$. (See \cite{monotonetriangles}.)
\item For $n \ge 1$, we have $\alpha(n;k_1,\ldots,k_n) = (-1)^{n-1} \alpha(n;k_2,\ldots,k_n,k_1-n)$. (A proof can be found in \cite{newproof}.)
\item For $n \ge 1$ and $p \ge 1$, we have 
$$
e_p(\Delta_{k_1},\ldots,\Delta_{k_n}) \alpha(n;k_1,\ldots,k_n) = 0. 
$$
(See Lemma~1 in \cite{newproof}.)
\end{enumerate}
The first property is obviously the analog of the shift-antisymmetry of ${L}_n({\mathcal T},{\bf k})$ as the latter can obviously be formulated as follows.
$$
(\id + E_{k_{i+1}} E^{-1}_{k_i} S_{k_i,k_{i+1}})  
{L}_n({\mathcal T},{\bf k}) = 0
$$
It is interesting to note that a special case of this property for 
$\alpha(n;k_1,\ldots,k_n)$
follows from Proposition~\ref{boundary}: if we specialize $k_{i+1}=k_i-1$ then
the first property simplifies to
\begin{multline*}
\alpha(n;k_1,\ldots,k_{i-1},k_i-1,k_i-1,k_{i+2},\ldots,k_n) 
+ \alpha(n;k_1,\ldots,k_{i-1},k_i,k_i,k_{i+2},\ldots,k_n) \\
- \alpha(n;k_1,\ldots,k_{i-1},k_i-1,k_i,k_{i+2},\ldots,k_n) = 0.
\end{multline*}
However, for integers $k_1, k_2, \ldots, k_{i},k_{i+2}, \ldots, k_n$ with 
$k_1 < k_2 < \ldots < k_{i-1} < k_i -1$ and 
$k_i < k_{i+2} < \ldots < k_{n-1} < k_n$, Proposition~\ref{boundary} implies this identity: in a monotone triangle $(a_{i,j})_{1 \le j \le i \le n}$ with bottom row $k_1,\ldots,k_{i-1},k_i-1,k_i,k_{i+2},\ldots,k_n$ we have either $a_{n-1,i}=k_i-1$, which corresponds to the case that we have $k_1,\ldots,k_{i-1},k_i-1,k_i-1,k_{i+2},\ldots,k_n$ as bottom row, or $a_{n-1,i}=k_i$, which corresponds to the case that $k_1,\ldots,k_{i-1},k_i,k_i,k_{i+2},\ldots,k_n$ is the bottom row.
As a polynomial in $k_1, k_2, \ldots, k_{i},k_{i+2}, \ldots, k_n$ is 
uniquely determined by its values on the set of these elements $(k_1,k_2,\ldots,k_{i},k_{i+2},\ldots,k_n) \in \mathbb{Z}^{n-1}$, the identity follows. 
 
\medskip

Concerning the second property, we have seen that it also holds for
${\mathcal L}_n({\mathcal T},{\bf k})$. Both properties  together actually 
imply \eqref{recursion}, see \cite{monotonetriangles}, and thus it would be interesting to give combinatorial proofs of these properties.

\subsubsection{Property (3) implies the refined alternating sign matrix theorem}
The third property is interesting as it holds also for 
Gelfand-Tsetlin patterns where it can easily be deduced from the shift-antisymmetry. However, it is 
a mystery that it also holds for monotone triangles, as we do not see how it can be deduced from the first property. Quite remarkably, it can be used to deduce the refined alternating sign matrix theorem as we explain next.

\medskip

The number $A_{n,i}$ of $n \times n$ alternating sign matrices, where the unique $1$ in the first row  is located in the 
$i$-th column is equal to the number of monotone triangles with bottom row $1,2,\ldots,n$ and 
$i$ appearances of $1$ in the first NE-diagonal, or, equivalently, the number of monotone triangles with 
bottom row $1,2,\ldots,n$ and $i$ appearances of $n$ in the last SE-diagonal. (This follows immediately from the standard bijection between alternating sign matrices and monotone triangles.)
If we assume that 
$k_1 \le k_2 < \ldots < k_n$, then the number of ``partial'' monotone triangles 
with $n$ rows, where the entries $a_{n,1}, a_{n-1,1}, \ldots, a_{n-i+1,1}$ are removed, no entry is smaller than $k_1$ and $a_{n,i} = k_i$ for $i=2,3,\ldots,n$ is equal to
$$
\left. (-1)^{i-1} \Delta_{k_1}^{i-1} \alpha(n;k_1,\dots,k_n)\right|_{(k_1,\ldots,k_n)=(1,1,2,\ldots,n-1)}.
$$
(A proof is given in \cite{refinedtrapezoids}.)
In fact, it follows quite easily by induction with respect to $i$ as
$$
- \Delta_{k_1} \left( \sum_{(l_1,\ldots,l_{n-1}) \in \mathbb{Z}^{n-1} \atop 
k_1 \le l_1 \le k_2 \le l_2 \le k_3 \le \ldots \le k_{n-1} \le l_{n-1} \le k_n, l_i \not= l_{i+1}} a(l_1,\ldots,l_{n-1})\right)  = 
 \sum_{(l_2,\ldots,l_{n-1}) \in \mathbb{Z}^{n-2} \atop 
 k_2 \le l_2 \le k_3 \le \ldots \le k_{n-1} \le l_{n-1} \le k_n, l_i \not= l_{i+1}} a(k_1,l_2,\ldots,l_{n-1}).
$$
This implies the first identity in 
$$
A_{n,i} = \left. (-1)^{i-1} \Delta_{k_1}^{i-1} \alpha(n;k_1,\dots,k_n) 
\right|_{(k_1,\ldots,k_n)=(1,1,2,\ldots,n-1)} = 
\left. \delta_{k_n}^{i-1} \alpha(n;k_1,\dots,k_n) 
\right|_{(k_1,\ldots,k_n)=(1,2,\ldots,n-1,n-1)}.
$$
The proof of the fact that the first expression is also equal to the last expression is similar. 
Therefore, by Property (3), 
\begin{multline*}
A_{n,i} = \left. (-1)^{i+n} \Delta_{k_1}^{i-1} \alpha(n;k_2,\dots,k_n,k_1-n) 
\right|_{(k_1,\ldots,k_n)=(1,1,2,\ldots,n-1)} \\ 
= \left. (-1)^{i+n} \delta_{k_1}^{i-1} E_{k_1}^{-2n+1+i} \alpha(n;k_2,\dots,k_n,k_1) 
\right|_{(k_2,\ldots,k_n,k_1)=(1,2,\ldots,n-1,n-1)}.
\end{multline*}
We use
$
E_x^{-m} = (\id - \delta_x)^m = \sum\limits_{j=0}^{m} \binom{m}{j} (-1)^{j} \delta_x^j
$
to see that this is equal to 
\begin{multline*}
\left. (-1)^{i+n} \delta_{k_1}^{i-1} \sum_{j=0}^{2n-1-i} \binom{2n-1-i}{j} (-1)^{j} \delta_{k_1}^j 
\alpha(n;k_2,\dots,k_n,k_1) \right|_{(k_2,\ldots,k_n,k_1)=(1,2,\ldots,n-1,n-1)} \\ =
\sum_{j=0}^{2n-1-i} \binom{2n-1-i}{j} (-1)^{i+j+n} A_{n,i+j}.
\end{multline*}
This shows that the refined alternating sign matrix numbers $A_{n,i}$ are a solution of the following system of linear equations.
$$
A_{n,i} = \sum_{k=1}^{n} \binom{2n-1-i}{k-i} (-1)^{k+n} A_{n,k}, \qquad 
1 \le i \le n
$$
In \cite{newproof}, it was shown that this system of linear equations together with the obvious symmetry $A_{n,i}=A_{n,n+1-i}$ determines the numbers $A_{n,i}$ inductively with respect to $n$. 

\medskip

It is worth mentioning that a similar reasoning can be applied to the doubly refined enumeration $\overline{\underline{A}}_{n,i,j}$ of $n \times n$ alternating sign matrices with respect to the position $i$ of the $1$ in first row and the position $j$ of the $1$ in the last row. This number is equal to the number of monotone triangles with bottom row $1,2,\ldots,n$ and $i$ appearances of $1$ in the first NE-diagonal and  $j$ appearances of $n$ in the last SE-diagonal, which implies (see 
\cite{refinedtrapezoids}) that
$$
\overline{\underline{A}}_{n,i,j} = 
\left. (-1)^{i-1} \Delta^{i-1}_{k_1} \delta^{j-1}_{k_n} \alpha(n;k_1,\ldots,k_n) \right|_{(k_1,\ldots,k_n)=(2,2,\ldots,n-1,n-1)}.
$$
Using the first and the third property of $\alpha(n;k_1,\ldots,k_n)$ displayed above we deduce 
the following identity.
\begin{multline*}
(\id + E^{n-1}_{k_n} E^{-n+1}_{k_1} S_{k_1,k_n}) V_{k_n,k_1} \alpha(n;k_1,\ldots,k_n) \\
= (-1)^{n-1} (\id + E^{n-1}_{k_n} E^{-n+1}_{k_1} S_{k_1,k_n}) V_{k_n,k_1} 
\alpha(n;k_2,\ldots,k_n,k_1-n) = 0
\end{multline*}
We apply $(-1)^{i-1} \Delta^{i-1}_{k_1} \delta^{j-1}_{k_n}$ to the equivalent identity
\begin{multline*}
0 = \alpha(n;k_1,\ldots,k_n) + \Delta_{k_1} \delta_{k_n} \alpha(n;k_1,\ldots,k_n) \\
+ E^{-2n+4}_{k_1} E^{2n-4}_{k_n} \alpha(n;k_n-n+3,k_2,\ldots,k_{n-1},k_1+n-3)
+ E^{-2n+4}_{k_1} E^{2n-4}_{k_n} \delta_{k_1} \Delta_{k_n} 
\alpha(n;k_n-n+3,k_2,\ldots,k_{n-1},k_1+n-3)
\end{multline*}
to see that
\begin{multline*}
0 = (-1)^{i-1} \Delta^{i-1}_{k_1} \delta^{j-1}_{k_n} \alpha(n;k_1,\ldots,k_n)
- (-1)^{i} \Delta^{i}_{k_1} \delta^{j}_{k_n} \alpha(n;k_1,\ldots,k_n) \\
+ E^{-2n+3+i}_{k_1} E^{2n-3-j}_{k_n}  (-1)^{i-1} \delta^{i-1}_{k_1} \Delta^{j-1} _{k_n}
\alpha(n;k_n-n+3,k_2,\ldots,k_{n-1},k_1+n-3) \\
+ E^{-2n+3+i}_{k_1} E^{2n-3-j}_{k_n}  (-1)^{i-1} \delta^{i}_{k_1} \Delta^{j} _{k_n}
\alpha(n;k_n-n+3,k_2,\ldots,k_{n-1},k_1+n-3).
\end{multline*}
Now we use the expansions
$$
E^{-2n+3+i}_{k_1} = (\id - \delta_{k_1})^{2n-3-i} = \sum\limits_{p=0}^{2n-3-i} \binom{2n-3-i}{p} (-1)^p \delta^p_{k_1} 
$$
and
$$
E^{2n-3-j}_{k_n} =  (\id + \Delta_{k_n})^{2n-3-j} = \sum\limits_{q=0}^{2n-3-j} \binom{2n-3-j}{q} \Delta^q_{k_n}
$$
to see that 
\begin{multline*}
0 = (-1)^{i-1} \Delta^{i-1}_{k_1} \delta^{j-1}_{k_n} \alpha(n;k_1,\ldots,k_n)
- (-1)^{i} \Delta^{i}_{k_1} \delta^{j}_{k_n} \alpha(n;k_1,\ldots,k_n) \\
+ \sum_{p=0}^{2n-3-i}  \sum_{q=0}^{2n-3-j} \binom{2n-3-i}{p} \binom{2n-3-j}{q} (-1)^{i-1+p} \delta^{p+i-1}_{k_1} \Delta^{q+j-1} _{k_n}
\alpha(n;k_n-n+3,k_2,\ldots,k_{n-1},k_1+n-3) \\
+  \sum_{p=0}^{2n-3-i}  \sum_{q=0}^{2n-3-j} \binom{2n-3-i}{p} \binom{2n-3-j}{q} (-1)^{i-1+p} \delta^{i+p}_{k_1} \Delta^{j+p} _{k_n}
\alpha(n;k_n-n+3,k_2,\ldots,k_{n-1},k_1+n-3).
\end{multline*}
We evaluate at $(k_1,k_2,\ldots,k_{n-1},k_n)=(2,2,3,\ldots,n-2,n-1,n-1)$ to arrive at 
$$
\overline{\underline{A}}_{n,i+1,j+1} - \overline{\underline{A}}_{n,i,j}  =
\sum_{p=0}^{2n-3-i}  \sum_{q=0}^{2n-3-j} \binom{2n-3-i}{p} \binom{2n-3-j}{q} (-1)^{i+j+p+q} 
\left( \overline{\underline{A}}_{n,q+j,p+i} - \overline{\underline{A}}_{n,q+j+1,p+i+1} \right).
$$
Computer experiments led us to the conjecture that this identity together with the obvious relations $\overline{\underline{A}}_{n,i,j} = \overline{\underline{A}}_{n,j,i}$ and $\overline{\underline{A}}_{n,i,j} = \overline{\underline{A}}_{n,n+1-i,n+1-j}$ determine  the doubly refined enumeration numbers $\overline{\underline{A}}_{n,i,j}$ uniquely inductively with respect to $n$.

\subsubsection{Property (1) and  (4) imply Property (3).} The analog of the fourth property is  true for Gelfand-Tsetlin tree sequences, see \eqref{toshow}, for which we gave a combinatorial proof in Section~\ref{diff}.
The significance of this property is that it can be used to deduce the third property from the first property. Since every symmetric polynomial in $X_1,X_2,\ldots,X_n$ can be written as a polynomial in the elementary symmetric functions, this property implies that
$$
p(E_{k_1},\ldots,E_{k_n}) \alpha(n;k_1,\ldots,k_n) = p(1,1,\ldots,1) 
\alpha(n;k_1,\ldots,k_n)
$$
for every symmetric polynomial $p(X_1,\ldots,X_n)$ in $X_1,\ldots,X_n$. This extends to symmetric polynomials in $X_1,X_1^{-1},\ldots,X_n,X^{-1}_n$: let $p(X_1,\ldots,X_n)$ be such a polynomial and $t \in \mathbb{Z}$ such that $p(X_1,\ldots,X_n) X^t_1 \cdots X^t_n=:q(X_1,\ldots,X_n)$ is a symmetric polynomial in $X_1,\ldots,X_n$ then
\begin{multline*}
p(E_{k_1},\ldots,E_{k_n}) \alpha(n;k_1,\ldots,k_n) = 
E^{t}_{k_1} \cdots E^{t}_{k_n} p(E_{k_1},\ldots,E_{k_n})  \alpha(n;k_1-t,\ldots,k_n-t) \\
= E^{t}_{k_1} \cdots E^{t}_{k_n} p(E_{k_1},\ldots,E_{k_n})  \alpha(n;k_1,\ldots,k_n) =  q(1,1,\ldots,1) \alpha(n;k_1,\ldots,k_n) = p(1,1,\ldots,1) \alpha(n;k_1,\ldots,k_n).
\end{multline*}
In particular, this shows that \eqref{toshow} is also true if all ``$\Delta$''s are replaced by ``$\delta$''s. Now we are ready to deduce Property (3) from Property (1) and Property (4): note that the operator $V_{x,y}$ is invertible as an operator on polynomials in $x$ and $y$: this follows as $V_{x,y}=\id + \delta_x \Delta_y$ and 
$$
V^{-1}_{x,y} = \sum_{i=0}^{\infty} (-1)^i \delta_x^i \Delta^i_y.
$$
(The sum is finite when applied to polynomials.)
Property (1) is obviously equivalent to 
$$
\alpha(n;k_1,\ldots,k_{i-1},k_{i+1}+1,k_i-1,k_{i+2},\ldots,k_n) = 
- V_{k_i,k_{i+1}} V^{-1}_{k_{i+1},k_i}  \alpha(n;k_1,\ldots,k_n).
$$
This implies 
\begin{multline*}
(-1)^{n-1} \alpha(n;k_2,\ldots,k_n,k_1-n) =
(-1)^{n-1} \alpha(n;k_2+1,\ldots,k_n+1,k_1-n+1) \\
= \prod_{i=2}^{n} V_{k_1,k_i} V^{-1}_{k_i,k_1} \alpha(n;k_1,\ldots,k_n).
\end{multline*}
Therefore, in order to show the third property, we have 
to prove that 
$$
\left( \prod_{i=2}^{n} V_{k_1,k_i}  - \prod_{i=2}^{n} V_{k_i,k_1} \right)
\alpha(n;k_1,\ldots,k_n) = 0.
$$
This follows from the fourth property as
\begin{multline*}
\prod_{i=2}^{n} V_{k_1,k_i}  - \prod_{i=2}^{n} V_{k_i,k_1} =
\prod_{i=2}^{n} (\id + \delta_{k_1} \Delta_{k_i}) -
\prod_{i=2}^{n} (\id + \Delta_{k_1} \delta_{k_i}) =
\sum_{r=0}^{n-1} \delta^r_{k_1} e_r(\Delta_{k_2},\ldots,\Delta_{k_n})
- \sum_{r=0}^{n-1} \Delta^r_{k_1}  e_r(\delta_{k_2},\ldots,\delta_{k_n}) \\
= \sum_{r=0}^{n-1} \left( \delta^r_{k_1} \left( e_r(\Delta_{k_1},\ldots,\Delta_{k_n}) -
\Delta_{k_1} e_{r-1}(\Delta_{k_2},\ldots,\Delta_{k_n}) \right) 
-  \Delta^r_{k_1} \left( e_r(\delta_{k_1},\ldots,\delta_{k_n}) -
\delta_{k_1} e_{r-1}(\delta_{k_2},\ldots,\delta_{k_n}) \right) \right) \\
= \sum_{r=0}^{n-1} \left(  \delta^r_{k_1}  e_r(\Delta_{k_1},\ldots,\Delta_{k_n})
- \Delta^r_{k_1}  e_r(\delta_{k_1},\ldots,\delta_{k_n}) \right) \\
- \sum_{r=1}^{n-1} \left(  \delta^r_{k_1} \Delta_{k_1} e_{r-1}(\Delta_{k_2},\ldots,\Delta_{k_n}) -  \Delta^r_{k_1} \delta_{k_1} e_{r-1}(\delta_{k_2},\ldots,\delta_{k_n}) \right) = \ldots \\
= \sum_{s=1}^{n} \sum_{r=1}^{n-s} (-1)^s \left( \Delta^{r+s-1}_{k_1} \delta^{s-1}_{k_1}
e_r(\delta_{k_1},\ldots,\delta_{k_n}) - \delta^{r+s-1}_{k_1} \Delta^{s-1}_{k_1}
e_r(\Delta_{k_1},\ldots,\Delta_{k_n}) \right).
\end{multline*}

\begin{appendix}

\section{The non-intersecting lattice paths point of view}
\label{path}

\begin{figure}
\begin{center}
\mbox{\scalebox{0.50}{%
    \includegraphics{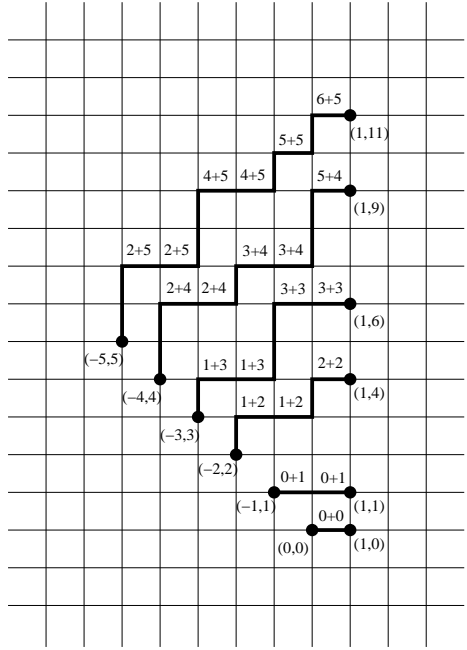}}}
\end{center}
\caption{\label{gesselfigure} Non-intersecting lattice paths}   
\end{figure}

In Figure~\ref{gesselfigure}, the family of non-intersecting lattice paths that corresponds to the Gelfand-Tsetlin pattern given in the introduction is displayed: in general, the lattice paths join the starting points $(0,0), (-1,1), \ldots, (-n+1,n-1)$ to
the end points $(1,k_1), (1,k_2+1), \ldots, (1,k_n+n-1)$, where the lattice paths can take east and north steps of length $1$ and end with a step to the east. 
As indicated in the drawing, the heights of the horizontal steps of the $i$-th 
path, counted from the bottom, can be obtained from the $i$-th southeast 
diagonal of the Gelfand-Tsetlin pattern, counted from the left, by adding $i$ to the entries in the respective diagonal of the Gelfand-Tsetlin pattern. By a well-known result on the enumeration of non-intersecting lattice paths of Lindstr\"om \cite[Lemma~1]{lindstroem} and of Gessel and Viennot 
\cite[Theorem~1]{gessel}, this number is equal to 
\begin{equation}
\label{det}
\det_{1 \le i, j \le n} \binom{k_j+j-1}{i-1},
\end{equation}
which is, by the Vandermonde determinant evaluation, equal to \eqref{product}.
(Note that $\binom{k_j+j-1}{i-1}$ is a polynomial in $k_j$ of degree $i-1$.)

\medskip

Interestingly, another possibility to extend the combinatorial interpretation of \eqref{product} to all $(k_1,\ldots,k_n) \in \mathbb{Z}_{\ge 0}^n$
is related to this interpretation in terms of families of non-intersection lattice paths: for arbitrary non-negative integers $k_1, k_2, \ldots, k_n$,
consider families of $n$ lattice paths with unit steps to the north and to the east (in general, these families are intersecting for the moment) that connect the starting points $(0,0), (-1,1), \ldots, (-n+1,n-1)$ to the endpoints $(0,k_1), (0,k_2+1), \ldots, (0,k_n+n-1)$, in any order. (Now we omit the vertical steps at the end of the paths.) Suppose that the $i$-th starting point $(-i+1,i-1)$ is connected to the $\pi_i$-th end point $(0,k_{\pi_i}+\pi_i-1)$ then the sign of the family is defined as the sign of the permutation 
$(\pi_1, \pi_2, \ldots, \pi_n)=\pi$. Then, \eqref{product} is the
signed enumeration of families of lattice paths with these starting points and end points. The merit of the theorem of Lindstr\"om and of Gessel and Viennot is the definition of a sign reversion involution on the families of {\it intersecting} lattice paths, which shows that only the non-intersecting families  remain in the signed enumeration. Depending on the relative positions of the numbers $k_1, k_2+1,\ldots,k_n+n-1$, there is only one permutation $\pi$ for which a family of non-intersecting lattice paths exists at all. This implies that the signed enumeration of families of lattice paths reduces essentially (i.e. up to the sign of $\pi$) to the plain enumeration of families of non-intersecting lattice paths. 

\medskip

\begin{figure}
\begin{center}
\mbox{\scalebox{0.50}{%
    \includegraphics{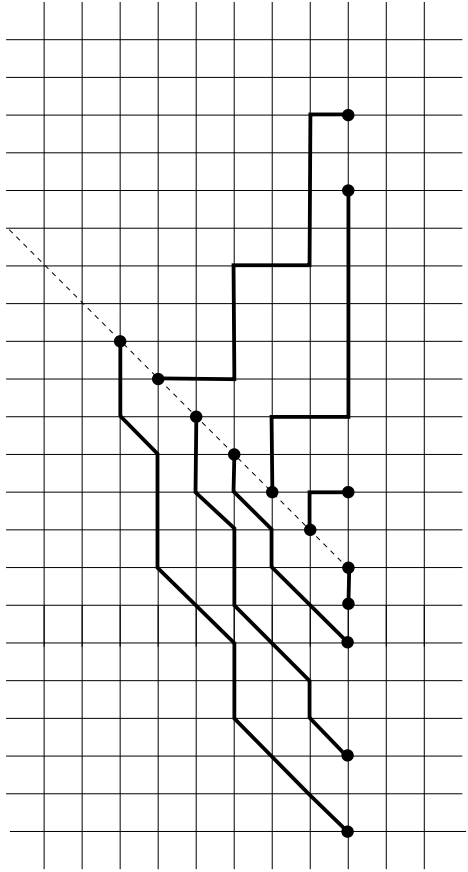}}}
\end{center}
\caption{\label{paths2} Non-intersecting lattice paths that go below the $x$-axis.}   
\end{figure}

Finally, it is worth mentioning (without proof) that the requirement that all $k_i$ are {non-negative} can be avoided. A close look at the proof shows that this requirement is useful at first place to guarantee that the location of the end points is not ``too far'' to the south of the starting points. If an end point is south-east of a starting point then there is obviously no lattice path connecting them which only uses steps of the form $(1,0)$ and $(0,1)$. However, in such a case it is convenient to allow steps of the form $(1,-1)$ and $(0,-1)$. Moreover if we require these paths to start with a step of the form $(0,-1)$ and let each step of the form $(1,-1)$ contribute a minus sign, we obtain an interpretation of \eqref{product} for all $(k_1,\ldots,k_n) \in \mathbb{Z}^n$. A typical situation is sketched in Figure~\ref{paths2}.

\section{Another proof of the shift-antisymmetry}
\label{shiftantisymmetry}

We sketch a 
(sort of) combinatorial proof of the shift-antisymmetry of the signed enumeration of 
Gelfand-Tsetlin {\it patterns}  with prescribed bottom row which does not rely  on
the notion of Gelfand-Tsetlin tree sequences. The argument is a bit involved and thus shows the 
merit of the notion of Gelfand-Tsetlin tree sequences. On the other hand,
it could be helpful for proving the analog property for monotone triangles as we have not established a notion 
that is analog to that of Gelfand-Tsetlin tree sequences for monotone triangles so far, 
see Section~\ref{concluding}.

\medskip

The following notion turns out to be extremely useful in order to avoid case distinctions: we define
$[x,y]:= \{z \in \mathbb{Z} | x \le z \le y\}$ if $x \le y$ as usual,  $[x,x-1] := \emptyset$ 
and  $[x,y] := [y+1,x-1]$  if $y+1 \le x-1$.  The latter situation is said to be an inversion. By considering all possible relative positions of $x, y, z$, it is not hard to see that 
$$
[x,y] \triangle [x,z+1] = [y+1,z+1],
$$
where $A \triangle B:= (A \setminus B) \cup (B \setminus A)$ is the symmetric difference.
In fact, concerning this symmetric difference, the following can 
be observed: either one set is contained in the other or the sets are disjoint. The latter situation 
occurs iff exactly one of $[x,y]$ and $[x,z+1]$ is an inversion. 
On the other hand, 
$$
[z,x] \triangle [y-1,x] = [x+1,z-1] \triangle [x+1,y-2] = [z,y-2] = [y-1,z-1]
$$
and we have $[z,x] \setminus [y-1,x] \not= \emptyset$ and $[y-1,x] \setminus [z,x] \not= \emptyset$ (which implies that the two sets are disjoint) iff exactly one of $[z,x]$ and $[y-1,x]$ is an inversion.

\medskip

Let 
${\mathcal L}_n(k_1,\ldots,k_n) := {\mathcal L}_n({\mathcal B},{\bf k})$ denote the set of Gelfand-Tsetlin 
patterns with bottom row $k_1,k_2,\ldots,k_n$ and ${L}_n(k_1,\ldots,k_n) := {L}_n({\mathcal B},{\bf k})$
the corresponding signed enumeration. The proof is by induction with respect to $n$. Nothing is to be done for $n=1$. Otherwise, it suffices to consider the case $j=i+1$. We fix $i \in \{1,2,\ldots,n-1\}$ and decompose ${\mathcal L}_n(k_1,\ldots,k_n)$ into four 
sets: let ${\mathcal L}^1_{n,i}(k_1,\ldots,k_n)$ denote the subset of patterns 
$(a_{p,q})_{1 \le q \le p \le n} \in {\mathcal L}_n(k_1,\ldots,k_n)$
 for which the replacement $a_{n,i} \to k_{i+1}+1$ and $a_{n,i+1} \to k_i-1$
produces another Gelfand-Tsetlin pattern (which is obviously an element of 
${\mathcal L}_n(k_1,\ldots,k_{i-1},k_{i+1}+1,k_i-1,k_{i+2},\ldots,k_n)$ then). If 
we perform this replacement we can either have a contradiction concerning the requirement for $l_{i-1}:=a_{n-1,i-1}$ or for $l_{i+1}:=a_{n-1,i+1}$. (There can not be a contradiction for 
$l_i:=a_{n-1,i}$ as $l_{i} \in [k_i,k_{i+1}]$ if and only if $l_{i} \in [k_{i+1}+1,k_i-1]$.) We let ${\mathcal L}^2_{n,i}(k_1,\ldots,k_n)$ denote the set of patterns, where we have a contradiction 
for $l_{i-1}$ but not for $l_{i+1}$,  ${\mathcal L}^3_{n,i}(k_1,\ldots,k_n)$ denote 
the set of patterns, where we have a contradiction for $l_{i+1}$ but not for $l_{i-1}$ and ${\mathcal L}^4_{n,i}(k_1,\ldots,k_n)$ denote the set of patterns, where we have a contradiction for both $l_{i-1}$ and $l_{i+1}$. Finally, we let ${L}^j_{n,i}(k_1,\ldots,k_n)$ denote the respective signed enumerations. We aim to show that 
\begin{equation}
\label{four}
{L}^j_{n,i}(k_1,\ldots,k_n) = -{L}^j_{n,i}(k_1,\ldots,k_{i-1},k_{i+1}+1,k_i-1,k_{i+2},\ldots k_n)
\end{equation}
if $j \in \{1,2,3,4\}$.

\medskip

The case $j=1$ is almost obvious, only the sign requires the following thoughts: 
having no contradiction for both $l_{i-1}$ and $l_{i+1}$ means that 
$l_{i-1} \in [k_{i-1},k_i] \cap [k_{i-1},k_{i+1}+1]$ and $l_{i+1} \in [k_{i+1},k_{i+2}] \cap [k_{i}-1,k_{i+2}]$. This is in fact true for patterns in ${\mathcal L}^1_{n,i}(k_1,\ldots,k_n)$ as well as for patterns in ${\mathcal L}^1_{n,i}(k_1,\ldots,k_{i+1}+1,k_i-1,\ldots,k_n)$. The intersection 
$[k_{i-1},k_i] \cap [k_{i-1},k_{i+1}+1]$ is empty if exactly one of the intervals is an inversion.
Thus we may assume that they are either both inversions or both not inversions. 
This implies that $l_{i-1}$ is an inversion for the patterns on the left if and only if it is an inversion for the patterns on the right. The same is true for $l_{i+1}$. On the other hand, 
$l_{i}$ is obviously an inversion on the left if and only if it is no inversion on the right, 
which takes care of the minus sign.

\medskip

We show \eqref{four} for $j=2$ (the case $j=3$ is analog by symmetry):
given an element of ${\mathcal L}^2_{n,i}(k_1,\ldots,k_n)$, we have  $l_{i-1}  \in [k_{i-1},k_i] \setminus [k_{i-1},k_{i+1}+1]$, whereas for an element of
${\mathcal L}^2_{n,i}(k_1,\ldots,k_{i+1}+1,k_i-1,\ldots,k_n)$, we 
have $l_{i-1}  \in [k_{i-1},k_{i+1}+1] \setminus [k_{i-1},k_i]$. The conditions for the other elements are the same. (In particular, $l_{i+1} \in [k_{i+1},k_{i+2}] \cap [k_{i}-1,k_{i+2}]$.) If we are in the case that either both sets $[k_{i-1},k_i]$ and $[k_{i-1},k_{i+1}+1]$ are no inversions or both sets are inversions then one set is contained in the other, which implies that one of the conditions for $l_{i-1}$ can not be met. However, then the condition for $l_{i-1}$ in the other set is that it lies in $[k_{i}+1,k_{i+1}+1]$. As the condition for $l_{i}$ is that it is contained in $[k_i,k_{i+1}]$ it follows, by the shift-antisymmetry for $n-1$, that the signed enumeration of the patterns in
this set must be zero. 

\medskip

If, however, exactly one set of $[k_{i-1},k_i]$ and $[k_{i-1},k_{i+1}+1]$ is an inversion then the sets are disjoint and their union is $[k_i+1,k_{i+1}+1]$. 
We decompose the two sets ${\mathcal L}^2_{n,i}(k_1,\ldots,k_n)$ and ${\mathcal L}^2_{n,i}(k_1,\ldots,k_{i+1}+1,k_i-1,\ldots,k_n)$ further according whether $l_i \in [k_{i-1}-1,k_i-1]$ or 
$l_{i} \in [k_{i-1}-1,k_{i+1}]$. (Observe that  also $[k_i,k_{i+1}]$ is the disjoint union of 
$[k_{i-1}-1,k_i-1]$ and $[k_{i-1}-1,k_{i+1}]$.) 
By the shift-antisymmetry for $n-1$, the signed enumeration 
of the elements in ${\mathcal L}^2_{n,i}(k_1,\ldots,k_n)$  which satisfy $l_i \in [k_{i-1}-1,k_i-1]$ is zero as the requirement for $l_{i-1}$ is that it is contained in $[k_{i-1},k_i]$.
Similarly, the signed enumeration of the elements in ${\mathcal L}^2_{n,i}(k_1,\ldots,k_{i+1}+1,k_i-1,\ldots,k_n)$ with $l_{i} \in [k_{i-1}-1,k_{i+1}]$ is zero. Thus, 
for the first set, we are left with the patterns that satisfy 
$l_{i-1} \in [k_{i-1},k_i]$ and  $l_{i} \in [k_{i-1}-1,k_{i+1}]$ and, for the second set,  
the patterns with $l_{i-1} \in [k_{i-1},k_{i+1}+1]$ and $l_i \in [k_{i-1}-1,k_i-1]$ remain. 
By the symmetry of these conditions and the shift-antisymmetry for $n-1$, we see that the signed enumeration of the first set 
is the negative of the signed enumeration of the second set: as for the sign observe that $l_{i-1}$ is an inversion on the left (which is the case iff $[k_{i-1},k_i]$ is an inversion) if and only if it is no inversion on the right (which is the case iff $[k_{i-1},k_{i+1}+1]$ is no inversion). The analog assertion is true for $l_{i}$ as it is an inversion on the left iff $[k_i,k_{i+1}]$ is an inversion and it is an inversion on the right iff $[k_{i+1}+1,k_i-1]$ is an inversion. Finally, for $l_{i+1}$ we have the situation that it is an inversion on the left iff it is an inversion on the right or the condition $l_{i+1} \in [k_{i+1},k_{i+2}] \cap [k_{i}-1,k_{i+2}]$ can not be met. 

\medskip

The case $j=4$ is similar though a bit more complicated and left to the interested reader.

\end{appendix}

\end{document}